\documentclass{amsart}[12]
\usepackage{a4wide} 
\usepackage{color}
\usepackage{hyperref}
\usepackage[foot]{amsaddr}
\usepackage{tensor}
\usepackage{enumitem} 
\usepackage{amsfonts,amssymb,amsmath,amsthm,graphicx,latexsym}
\usepackage{mathrsfs,dsfont}
\usepackage{pgfplots}
\usepackage{tikz}
\usepackage{bm}
\usepackage{caption}
\newcommand{\A}{\mathcal{A}}
\newcommand{\R}{\mathbb R}
\newcommand{\C}{\mathbb C}
\newcommand{\N}{\mathbb N}
\newcommand{\Z}{\mathbb Z}
\newcommand{\M}{\mathbb M}
\newcommand{\vo}{\mathtt{v}}
\newcommand{\vol}{\mathtt{vol}}
\newcommand{\vp}{\varphi}
\newtheorem{thm}{Theorem}[section]
\newtheorem{prop}[thm]{Proposition}
\newtheorem{lemm}[thm]{Lemma}
\newtheorem{cor}[thm]{Corollary}
\newtheorem{rem}[thm]{Remark}
\newtheorem{exmp}[thm]{Examples}
\newtheorem{exm}[thm]{Example}
\newtheorem{df}[thm]{Definition}
\newcommand{\ve}{\varepsilon}

\newcommand{\Symz}{{\mathtt{Symm}}_d^+(\Z)}
\newcommand{\Sym}{{\mathtt{Symm}}_d^+(\R)}
\newcommand{\Symm}{{\mathtt{Symm}}_2^+(\R)}

\newcommand{\NC}{\mathcal{NC}}
\newcommand{\Pa}{\mathcal{P}}

\newcommand{\bmo}{\mathrm{bmo}}

\newcommand{\I}{\mathbf {I}}
\newcommand{\D}{\mathtt {D}}
\newcommand{\X}{\mathbb{X}}
\newcommand{\XX}{\mathtt{X}}
\newcommand{\B}{\mathtt {B}}
\newcommand{\BB}{\mathtt {BM1}}
\newcommand{\BC}{\mathtt {BM2}}

\newcommand{\bproof}{\noindent{\bf Proof: }}
\newcommand{\eproof}{\hfill $\Box$\\}

\newcommand{\tp}{\mathbin{\hbox{$\bigcirc$\hbox to
      0pt{\hspace{-0.81em}$\scriptstyle\top$\hfil}}}}
\begin{document}
\title[bm-Poisson type limit theorems]%
  {Analogues of Poisson type limit theorems in discrete bm-Fock spaces}

\author[L. Oussi]{Lahcen Oussi $^{1,2}$}
\address{$^1$ Institute of Mathematics, Wroclaw University, Pl. Grunwaldzki 2, 50-384 Wroclaw, Poland}
\address{$^2$lahcen.oussi@math.uni.wroc.pl}
\address{$^3$jwys@math.uni.wroc.pl}
\author[J. Wysoczański]{Janusz Wysoczański $^{1,3}$}
\begin{abstract}
We present analogues of the Poisson limit distribution for the noncommutative bm-independence, which is associated with several positive symmetric cones. We construct related discrete Fock spaces with creation, annihilation and conservation operators, and prove Poisson type limit theorems for them. Properties of the positive cones, in particular the \textit{volume characteristic} property they enjoy, and the combinatorics of labelled noncrossing partitions, play crucial role in these considerations. 
\end{abstract}


\maketitle
\section{Introduction}
In this paper we study analogues of the classical Poisson limit theorem in the noncommutative probability framework set by bm-independence. Our study can be considered as parallel to the results of Crismale, Griseta and the second-named author \cite{CGW2020} and \cite{CGW2021}, where the weakly monotone case has been considered. Here we develop theory extending the monotone case of Muraki \cite{Mur0}. 

In general, a noncommutative probability space is a unital $*-$algebra $\A$ (or a $C^*-$algebra) with a normalized state $\vp$ on it, which plays the role of classical expectation. 

There are several notions of independence in such framework, which are the rules for computing mixed moments, i.e. expressions of the form $\vp(a_{i_1}a_{i_2}\ldots a_{i_k})$ for $a_{i_1}, a_{i_2}, \ldots , a_{i_k}\in \A$. The universal ones are the freeness of Voiculescu \cite{Voi}, the monotonic independence of Muraki \cite{Mur1} and the Boolean independence \cite{S.W}. There are other notions of noncommutative independence, which are mixtures of these. For more information on this topic we refer to the Introduction in \cite{bmPoisson}. 

Here we develop the theory of bm-independence introduced in \cite{JW2008, JW2009}, which combines Boolean and monotonic ones. A particular  feature of this notion is that it is defined for random variables indexed by a partially ordered set. The definition of this notion is presented in the next section.   In this study, however, we consider only specific partially ordered sets, which are defined by \textit{positive symmetric cones} in Euclidean spaces, classification of which can be found in \cite{J.F.A.K}. These cones have an additional geometric property called the \textit{volume characteristic}, which allows to develop analogues of classical properties in such general framework. 

Noncommutative analogues of the classical Poisson limit theorem arise in two ways. The first one is by considering the convolution powers of Bernoulli type distributions, and this can be done for the free, monotone and boolean convolutions. This can be generalized to arrays of independent random variables (see \cite{Sp} for free independence and \cite{bmPoisson} for bm-independence). The second method, used by Muraki \cite{Mur0}, comes from constructions of creation $A_i^+$, annihilation $A_i^-$ and conservation $A_i^{\circ}$ operators ($i\in\N$) on appropriate Fock space associated, as the toy model, with a given noncommutative independence. In particular, these operators are independent (in the given noncommutative sense). For a parameter $\lambda \geq 0$, one considers the sums 
\begin{equation}\label{Poisson general}
 S_{N}(\lambda):=\frac{1}{\sqrt{N}}\sum_{i=1}^{N}(A_{i}^{+}+A_{i}^{-})+\lambda\sum_{i=1}^{N}A_{i}^{\circ}
\end{equation}
and one studies the limits of moments (for every $p\in\N$)
\[
m_p(\lambda):=\lim_{N\to +\infty} \vp((S_N(\lambda))^p). 
\]
Then these are moments of Poisson type measure (which is to be identified in each case) associated with the given noncommutative independence. This type of Poisson limit is motivated by the Fock space of the classical Poisson process \cite{Mey}. In particular, for free independence this is the free Poisson (i.e. Marchenko-Pastur) measure \cite{Sp}, whereas for monotonic independence the measure has been identified by Muraki \cite{Mur0}. Our constructions and results extend those  of Muraki, so we recall his description of moments of the limit measure in the monotonic case (for details we refer to \cite{Mur0}):
\begin{equation}\label{convergence}
m_{p}(\lambda):=\lim_{N\rightarrow +\infty}\vp\left((S_{N}(\lambda))^{p}\right)=\sum_{\pi \in \NC_{2}^{1, i}(p)}V(\pi)\lambda^{s(\pi)},
\end{equation}
where $\NC_{2}^{1, i}(p)$ denotes the set of all noncrossing  partitions of the set $[p]:=\{1, \ldots, p\}$, with pair blocks and inner singleton blocks. In addition, $V(\pi)$ counts (approximately, as $N\rightarrow +\infty$) the ratio of all  blocks of $\pi\in \NC_{2}^{1, i}(p)$ with monotonic labellings by integers $\{1, 2, \ldots, N\}$,  divided by $N^{\frac{p-s(\pi)}{2}}$, and $s(\pi)$ is the number of singleton blocks in $\pi$. 

The independence we study here can be defined for any partially ordered set, however such generality does not allow for specific results. In particular one needs additional structure of the poset. One of the most interesting ones is guaranteed by partial orders defined by positive symmetric cones, classification and properties of which can be found in \cite{J.F.A.K}. 

More in detail we consider the following three classes of positive symmetric cones $\Pi_d$ (with $d\in\N$) in particular vector spaces $\XX_d$:
\begin{enumerate}
  \item The cones $\Pi_d={\R}_{+}^{d}:=\{ \xi=(\xi_{1}, \xi_{2}, \ldots, \xi_{d})\in {\R}^{d}: \xi_{1}, \ldots, \xi_{d}\geq 0\}$ in the vector spaces ${\XX_d}={\R}^{d}$, where for $\xi=(\xi_{1}, \ldots, \xi_{d}), \rho=(\rho_{1},\ldots, \rho_{d})\in {\R}^{d}$, the partial order $\preceq$ is given by coordinatewise comparison:
$$\xi\preceq \rho   \quad \text{if}\quad  \xi_{1}\leq \rho_{1}, \ldots , \xi_{d}\leq \rho_{d}.$$
  \item The Lorentz light cones $\Lambda_{d}^{1}$ in $(d+1)$-Minkowski space time ${\XX_d}={\R}_{+}\times{\R}^{d}$ with the positive cone $\Pi_d=\Lambda_{d}^{1}:=\left\{(t; x_{1}, ..., x_{d})\in {\X_d}: t\geq \left(\sum_{i=1}^{d}x_{i}^{2}\right)^{\frac{1}{2}}\right\}$. In this case, for $\xi= (t; x_{1}, ..., x_{d})$, $\rho=(s; y_{1}, ..., y_{d})\in {\R}_{+}\times{\R}^{d}$, we have
$$\rho\preceq \xi   \quad \text{if}\quad    t-s\geq \left(\sum_{i=1}^{d}(x_{i}-y_{i})^{2}\right)^{\frac{1}{2}}.$$
   \item The cones $\Pi_d=\Sym$ of positive semidefinite real symmetric $(d\times d)$ matrices:
\[
\Sym:=\{\xi \in \M_{d}({\R}): \xi =\xi^T,\ \xi\geq 0 \} \subset \XX_d=\M_{d}({\R}).
\]
For $\xi, \rho\in\M_{d}({\R})$, one defines 
\[
\rho \preceq \xi  \quad \text{ if } \quad \xi-\rho\in\Sym.
\]  
\end{enumerate}
These cones satisfy the property of \emph{volume characteristic} which plays crucial role in our studies. This property, introduced in \cite{J.A}, reads as follows:
\begin{thm}[volume characteristic \cite{J.A}]\label{volchar}
 For each of the positive symmetric cones $\Pi_d$ we consider, there exists a sequence $(\gamma_{m}(\Pi_d))_{m\geq 1}$ such that for any $\xi\in \Pi_d$ and any $m\in{\N}$
 $$\gamma_{m}(\Pi_d)=\frac{1}{\vo(\xi)^{m}}\int_{\rho\in [0, \xi]}\vo(\rho)^{m-1}d(\rho),$$
 where $\vo(\xi):=\vol[0, \xi]$ denotes the Euclidean volume of the interval $[0, \xi]\subset \Pi_d\subset {\R}^{m}$ and $d(\rho)$ is the Lebesgue measure on ${\R}^{m}$ (the dimension $m$ is minimal for the embedding $\Pi_d\subset {\R}^{m}$).
\end{thm}
The \textit{volume characteristic} property generalizes the following simple fact for the positive cone ${\R}_{+}$, where $\vo(x):=\vol[0, x]=x$ for $x>0$:
\[
\frac{1}{\vo(x)^m} \int_{y\in[0,x]} \vo(y)^{m-1} dy   = \frac{1}{x^m}\int_{0}^{x}y^{m-1}dy=\frac{1}{m}=\gamma_{m}({\R}_{+}),
\]
In what follows we will write simply $\gamma_m$ if the cone is specified. 

In our studies of Poisson type limit theorems we deal with bm-independent random variables, indexed by elements of positive symmetric cones $\Pi_d$. This creates some challenges regarding the formulation of the theorems. These challenges are the same as those given in our study \cite{bmPoisson} of the law of small numbers for bm-independent random variables. For the reader's convenience we recall briefly how these problems has been solved in \cite{bmPoisson}.  
\begin{enumerate}
\item Instead of the index set of positive integers we use the following discrete subsets ${\I}\subset\Pi_d$ depending on the positive cone $\Pi_d$:
\begin{enumerate}
 \item ${\I}:={\N}^{d}$ if $\Pi_d={\R}_{+}^{d}$,
 \item ${\I}:={\N}\times {\Z}^d$ if $\Pi_d=\Lambda_{d}^{1}$ is the Lorentz light cone,
 \item ${\I}:=\Symz$ if $\Pi_d=\Sym\subset{\M}_{d}({\R})$.
\end{enumerate}
\item The range of summation $j\in \{1,2,\ldots , N\}=[1,N]\cap \N$ in the partial sum $\sum\limits_{j=0}^{N}A_{j}^{\ve}$ (with $\ve\in \{-, +, \circ\}$) in the formula \eqref{Poisson general}, is replaced by the finite summation $\sum\limits_{\xi \in [0, \rho]_{{\I}}}A_{\xi}^{\ve}$, where $[0, \rho]_{{\I}}:=[0, \rho]\cap {\I}$ for $ \rho\in\Pi_d $ and $[0, \rho]:=\{\xi\in \Pi_d: 0\preceq \xi\preceq \rho\}$ is an interval in the positive  cone $\Pi_d$.
\item The replacement of the convergence $N\rightarrow +\infty$ in the formula \eqref{convergence} has to be formulated separately for each positive symmetric cone under consideration. We use the notation $\rho\xrightarrow[]{\Pi}\infty$ in the following cases:
\begin{itemize}
   \item  For $\Pi_d={\R}_{+}^{d}$ and $\rho:=(a_{1}, \cdots ,a_{d})\in \Pi_d$, then $\rho \xrightarrow[]{\Pi}\infty$ means that $a_{j}\rightarrow \infty$ for all $1\leq j\leq d$,
\item For $\Pi_d=\Lambda_{d}^{1}$ and $\rho:=(t;x_{1}, \ldots, x_{d})\in \Pi_d$, then $\rho\xrightarrow[]{\Pi}\infty$  means that  $t-\left(\sum_{i=1}^{d}x_{i}^{2}\right)^{\frac{1}{2}}\rightarrow \infty$,
\item For $\Pi_d=\Sym$ and $\rho \in \Pi_d$, if $0<\rho_{1}\leq \rho_{2}\leq \cdots \leq \rho_{d}$ are the eigenvalues of $\rho$, then $\rho\xrightarrow[]{\Pi}\infty$ means that  $\rho_{1}\rightarrow \infty$ (and consequently $\rho_{j}\rightarrow \infty$ for all $1\leq j \leq d$).
  \end{itemize}

\item
For a function $f: \Pi_d\mapsto {\R}$, we define $\lim\limits_{\rho \xrightarrow[]{\Pi}\infty}f(\rho)=\alpha$ if for each $\epsilon>0$ there exists $\mu\in\Pi_d$ such that for every $\rho\in\Pi_d$ if $\mu\preceq \rho$ then $|f(\rho)-\alpha|<\epsilon$.
\item Furthermore, the normalization factor $\sqrt{N}$ in \eqref{Poisson general} is replaced by $\sqrt{\#[0, \rho]_{\I}}$, where $\#[0, \rho]_{\I}$ is the number of elements in $[0, \rho]_{\I}:=[0, \rho]\cap \I$, which is asymptotically the same as the Euclidean volume $\vo(\rho):=\vol[0, \rho]$ of an interval $[0, \rho]\subset\Pi_d$. 

For the reader's convenience we recall the following formulas for the Euclidean volume of intervals $[0, \rho]$ (cf.\  \cite{J.F.A.K}):
\begin{enumerate}
\item $\displaystyle \vo(\rho)=\prod_{j=1}^{d}a_j$ if $\rho:=(a_1, \ldots , a_d)\in\Pi_d=\R_+^d$,
\item $\displaystyle \vo(\rho)=\alpha_d(t^2-\|x\|^2)^{\frac{d+1}{2}}$ for some constant $\alpha_d$, if $\rho:=(t; x) \in \Pi_d = \Lambda_d^1$,
\item $\displaystyle \vo(\rho)=\beta_d\left(\prod_{j=1}^{d}\lambda_j\right)^{\frac{d+1}{2}}=\beta_d\Big(\det(\rho)\Big)^{\frac{d+1}{2}}$  for some constant $\beta_d$, if $(\lambda_1, \ldots , \lambda_d)$ are the eigenvalues of $\rho\in \Pi_d={\Sym}\subset {\M}_{d}(\R)$. 
\end{enumerate}
\end{enumerate}
With these changes our bm-analogue of the formula \eqref{Poisson general} can be written as
\begin{equation}\label{bmP}
S_{\rho}(\lambda):=\frac{1}{\sqrt{\vo(\rho)}}\sum_{\xi\in[0, \rho]_{\I}}(A_{\xi}^{+}+A_{\xi}^{-})+\lambda\sum_{\xi\in[0, \rho]_{\I}}A_{\xi}^{\circ}.
\end{equation}

The main goal of our study is the limit 
\[
\lim_{\rho\xrightarrow[]{\Pi}\infty}\vp((S_{\rho}(\lambda))^p), \quad p\in{\N},
\]
which we shall describe in a combinatorial manner. 

Our starting point is the construction of \emph{discrete bm-Fock spaces} (for each positive symmetric cone) and the creation $A_{\xi}^{+}$, annihilation $A_{\xi}^{-}$ and conservation $A_{\xi}^{\circ}$ operators on it. The construction  is related to the bm-product Fock space \cite{JW1} and  generalizes the ideas of Muraki. 

In \cite{J.W3} the second-named author provided a model for bm-independence, where \textit{bm-Fock space} has been constructed and the algebras $\{\A_{\xi}: \xi\in \X\}$, generated by creation $A^+_{\xi}:=A^+(e_{\xi})$ and annihilation $A^-_{\xi}:=A^-(e_{\xi})$ operators, are bm-independent if the vectors $\{e_{\xi}:\xi\in \X\}$ are mutually orthogonal.  This construction, however, is not suitable for the present study, which requires a discrete version of it. Therefore, in the next section, we present new constructions  of \textit{discrete bm-Fock spaces}, each of which is related with the given positive symmetric cone. 
\section{Preliminaries}
In this section we present the framework  for our study.

A \emph{noncommutative probability space} $(\mathcal{A}, \vp)$  consists of a unital $C^*$-algebra $\mathcal{A}$ and a state $\vp$ on it, that is, a linear functional $\vp: \mathcal{A}\rightarrow {\C}$ which is positive $(\vp(a^*a)\geq 0$ for all $a\in\mathcal{A})$ and unital $(\vp(1_{\mathcal{A}})=1$ for the unit $1_{\mathcal{A}}\in\mathcal{A})$. Self-adjoint elements $a=a^{*}\in \mathcal{A}$ are called \emph{noncommutative random variables} and the distribution  of a random variable $a=a^*\in  \A$ is a probability measure $\mu$ with the moments given by the sequence $\displaystyle m_n(\mu):=\vp(a^n)$. It always exists, since the sequence of moments is positive definite. 

General references on noncommutative probability are \cite{HO, Mey, Par} and references therein.
\subsection{bm-independence}
The general formulation of \emph{bm-independence} was given by the second-named author in \cite{J.W3} for families of algebras indexed by partially ordered sets. If $(\X, \preceq)$ is a poset, we shall use the following notation: $x\sim y$ if $x,y\in \X$ are comparable, $x\nsim y$ if $x,y\in \X$ are incomparable and $x\prec y$ if $x\preceq y$ and $x\neq y$.
\begin{df}[bm-independence]
Let $(\mathcal{A}, \vp)$ be a noncommutative probability space. We say that a family $\{\mathcal{A}_{\xi}: \xi\in {\X}\}$ of subalgebras of $\mathcal{A}$, indexed by a partially ordered set $({\X}, \preceq)$, is \textbf{bm-independent} in $(\mathcal{A}, \vp)$ if the following conditions hold: 
\begin{itemize}
 \item ${\BB}:$ If $\xi, \rho, \eta \in {\X}$ satisfy: $\xi\prec \rho\succ \eta$ or $\xi\nsim \rho\succ  \eta$ or $\xi \prec \rho \nsim\eta$, then for any $a_{1}\in \mathcal{A}_{\xi}, a_{2}\in \mathcal{A}_{\rho}, a_{3}\in \mathcal{A}_{\eta}$ we have 
 \begin{equation}\label{bm1d}
  a_{1}a_{2}a_{3}=\vp(a_{2})a_{1}a_{3}.
 \end{equation}
 \item ${\BC}:$ If $\xi_{1}\succ \cdots \succ \xi_{m}\nsim \cdots \nsim \xi_{k}\prec \xi_{k+1}\prec \cdots \prec \xi_{n}$ for some $1\leq m\leq k\leq n$ and $\xi_{1}, \cdots, \xi_{n}\in {\X}$, with $a_{j}\in \mathcal{A}_{\xi_{j}}$ for $1\leq j\leq n$, then 
 \begin{equation}\label{bm2d}
  \vp(a_{1} \cdots a_{n})=\prod_{j=1}^{n}\vp(a_{j}).
 \end{equation}
\end{itemize}
\end{df}
Noncommutative random variables $\{a_{\xi}\in \mathcal{A}: \xi\in {\X}\}$ are called \emph{bm-independent} if the subalgebras $\mathcal{A}_{\xi}$ they generate are bm-independent. For more information about properties of bm-independence we refer to \cite{J.W3, J.A}.
\begin{rem}
In the above definition, if ${\X}$ is totally ordered (i.e. every two elements in ${\X}$ are comparable), then we obtain the monotone independence. On  the other hand, if $\X$ is totally disordered (i.e. none of the elements of ${\X}$ are comparable), then we obtain the Boolean  independence.\\
The two conditions ${\BB}$ and ${\BC}$ allow to compute all joint moments $\vp(a_{1}\cdots a_{n})$ for bm-independent random variables $a_{1}, \ldots, a_{n}$ via marginals $\vp_{j}:=\vp|\mathcal{A}_{j}$, i.e. by the restriction of $\vp$ to the subalgebras they generate (c.f. \cite{J.W3}, Lemmas 2.3, 2.4). 
The idea is that first one applies ${\BB}$ as many times as possible, and then what remains is subject to ${\BC}$. We also refer to \cite{bmPoisson} for an algorithm which allows to evaluate explicitly  joint moments using the conditions ${\BB}$ and ${\BC}$.
\end{rem}
\subsection{Noncrossing partitions with pair or singleton blocks and bm-orders}
For $p\in {\N}$, a \emph{partition} $\pi$ of the finite set  $[p]:=\{1, \ldots, p\}$ is a collection of disjoint  nonempty subsets of $[p]$, called \emph{blocks}, whose union is $[p]$. For a partition $\pi$ with $k$ blocks $B_{1}, \ldots, B_{k}$ we will write  
$\pi:=(B_{1}, \ldots,  B_{k})$ to indicate that the blocks are ordered by their minimal elements: $\min(B_j)<\min(B_{j+1})$ for $1\leq j \leq k-1$. We denote by $\Pa(p)$ the set of all partitions $\pi$ of $[p]$, and by $|B|$ the cardinality of a block $B$ in $\pi$. If $|B|=1$, then the block $B$ consists of one element and is called \emph{singleton}. The numbers of blocks in $\pi$ will be denoted by $b(\pi)$ and the number of  singletons of $\pi$ by $s(\pi)$. One says that a partition $\pi$ has a \emph{crossing} if there exists two distinct blocks $B_{i}$ and $B_{j}$ in $\pi$, and elements $u_{1}, u_{2}\in B_{i}, v_{1}, v_{2}\in B_{j}$ such that $u_{1}< v_{1}<u_{2}<v_{2}$. Otherwise, $\pi$ has \emph{no crossings} and it is called a \textit{noncrossing partition}. The set of all \emph{noncrossing partitions} of $[p]$ will be denoted by $\NC(p)$ (or by $\NC(p, k)$ for noncrocssing partitions in $\NC(p)$ with exactly $k$ blocks).\\
Pictorially, a partition $\pi\in\Pa(p)$ will be presented by drawing the integers 1 to $p$ (from right to left) on a horizontal line and then joining the elements of a blocks by lines above; in particular singletons are marked by vertical lines. For instance, the graphical representation of the partition $$\pi=(\{1\}, \{2, 4\}, \{3, 5\}, \{6, 11\}, \{7, 9, 10\}, \{8\})\in \Pa(11),$$ is the following figure:
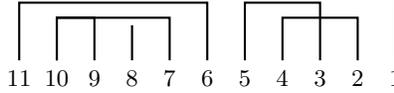
\begin{center}
 \begin{tikzpicture}[thick,font=\small]
 \path (0,0) node[] (a) {11}
  (0.5,0) node[] (b) {10}
 
            (1,0) node[] (c) {9}
 
            (1.5,0) node[] (d) {8}
 
            (2,0) node[] (e) {7}
 
            (2.5,0) node[] (f) {6}
            (3,0) node[] (g) {5}
 
            (3.5,0) node[] (h) {4}
 
            (4,0) node[] (i) {3}
 
            (4.5,0) node[] (j) {2}
 
            (5,0) node[] (k) {1}; 
      \draw (a) -- +(0,1) -| (f);
      \draw (g) -- +(0,1) -| (i);
       
       \draw (k)-- +(0,1);
        \draw (d)-- +(0,0.7);
     
   \draw (b) -- +(0,0.8) -| (e);
   \draw (b) -- +(0,0.8) -| (c);
    \draw (h) -- +(0,0.8) -| (j);

    \end{tikzpicture}  
 \end{center}
This partition has two singletons $\{1\}$ and $\{8\}$ and a crossing between the blocks $\{2,4\}$ and $\{3,5\}$. Of course, singletons do not produce crossings.
\begin{df}[Inner and outer blocks of a noncrossing partition]
A block $B_{j}\in\pi\in\NC(p)$ is called \textbf{inner} if there exists another block $B_{i}$ such that \begin{equation*}
\min B_{i}< \min B_{j}\leq \max B_{j}< \max B_{i},
\end{equation*}
where $\min B$ (resp. $\max B$) denotes the minimal (resp. maximal) element of the block $B$.
Otherwise, $B_j$ is called an \textbf{outer} block.
\end{df}
Observe that in the noncrossing partition every block is either inner or outer. 

It is convenient to define a partial order $\preceq_{\pi}$ on the blocks of a noncrossing partition $\pi=(B_{1}, \ldots, B_{k})\in \NC(p, k)$. Namely, we will write $B_{i}\preceq  _{\pi} B_{j}$  if $B_{j}$ is inside  $B_{i}$, that is  
\begin{equation*}
\min B_{i}< \min B_{j}\leq \max B_{j}< \max B_{i}.
\end{equation*}
A noncrossing partition $\pi=(B_{1}, \ldots, B_{k})\in \NC(p, k)$ is called a \emph{pair partition} (notation: $\pi \in \NC_{2}(p)$) if $|B_{j}|=2$ for all $1\leq j\leq k$ (each block $B_{i}$ contains exactly  two elements). On the other hand, a partition $\pi\in \NC(p)$ is called a \emph{noncrossing partition with pair or  singleton blocks} if $ |B|\in \{1, 2\}$ for all $B\in\pi$; the set of all such partitions will be denoted by $\NC_{2}^{1}(p)$. We distinguish two subsets of $\NC_{2}^{1}(p)$ which will appear in the combinatorial descriptions of our Poisson type limit  distribution in the next section:
\begin{itemize}
 \item  The set $\NC_{2. o}^{1, i}(p)\subset \NC_{2}^{1}(p)$, in which the partitions have \textbf{no inner pair blocks} and \textbf{no outer singletons}, i.e. pair blocks must be outer and singletons must be inner.
 \item The set $\NC_{2}^{1, i}(p)\subset \NC_{2}^{1}(p)$ in which the partitions have \textbf{no outer singletons} (i.e. singletons must be inner blocks, but pair blocks can be either inner or outer).
 \end{itemize}
Let $(\xi_{p}, \ldots, \xi_{1})$ be a sequence of elements of ${\X}$. For $\xi\in\{\xi_{p}, \ldots, \xi_{1}\}$, we define a subset $B(\xi):=\{1\leq j\leq p: \xi_{j}=\xi\}\subset [p]$, and denote by $i=i(\xi):=\min B(\xi)$ being the minimal element. If the cardinality of the set $\{\xi_{p}, \ldots, \xi_{1}\}$ is $k$, then we obtain $k$ disjoint subsets $B_1, \ldots , B_k$, which form a partition of $[p]$. Re-ordering the blocks to get that $\min(B_i)<\min(B_{i+1})$ for $1\leq i \leq k-1$, we obtain then a partition $\pi=(B_{1}, \ldots, B_{k})\in \Pa(p)$. In this case, we say that the partition $\pi$  is \textbf{adapted} to the  sequence $(\xi_{p}, \ldots, \xi_{1})$ and we denote this by $(\xi_{p}, \ldots, \xi_{1})\sim\pi$. Observe that the partition adapted to a sequence is unique. 

For a finite subset $\mathbf{J}\subset\I$ and a partition $\pi=(B_{1}, \ldots, B_{k})\in\Pa(p)$ with $k$ blocks, we define the \textbf{label function} $\mathtt{L}: \pi\longrightarrow \mathbf{J}$ by $\mathtt{L}(B_{i})\in\mathbf{J}$ for $1\leq i\leq k$. If we consider $(\xi_{p}, \ldots, \xi_{1})\sim\pi$ as above, then for $\xi\in\{\xi_{p}, \ldots, \xi_{1}\}$, $\mathtt{L}(B(\xi)):=\xi$ will be called \textbf{the label} of the block $B(\xi)$, and $(\mathtt{L}(B_{k}), \ldots, \mathtt{L}(B_{1}))$ will be called \textbf{the label sequence} of the partition $\pi$.
\begin{rem}\label{equalsets}
If a partition $\pi=(B_{1}, \ldots, B_{k})\in \Pa(p)$  is adapted to  a sequence $(\xi_{p}, \ldots, \xi_{1})$, then  the set $\{\xi_{p}, \ldots, \xi_{1}\}$ of elements and the set $\{\mathtt{L}(B_{k}), \ldots, \mathtt{L}(B_{1})\}$ of labels are equal. Moreover, the label sequence $(\mathtt{L}(B_{k}), \ldots, \mathtt{L}(B_{1}))$ is uniquely defined.
\end{rem}
Following \cite{CGW2021}, let $\pi=(B_{1}, \ldots, B_{k})\in\NC_{2}^{1}(p)$ be a given partition and let $B_{i}$ and  $B_{j}$ be two blocks, which for any other block $B_{l}\in\pi$, satisfy
$$\text{if $B_{i}\preceq_{\pi}B_{l}\preceq_{\pi} B_{j}$ then $B_{i}=B_{l}$ or $B_{j}=B_{l}$}.$$
Then we say that the block $B_{j}$ is a $\textbf{direct successor}$ of the block $B_{i}$; equivalently, the block $B_{i}$ is a $\textbf{direct predecessor}$ of the block $B_{j}$.
\begin{exm}
Consider $\pi=(B_{1}, B_{2}, B_{3}, B_{4}, B_{5})\in\NC_{2}(8)$ with $B_{1}=\{1, 8\}, B_{2}=\{2, 4\}, B_{3}=\{3\}, B_{4}=\{5, 7\}$ and $B_{5}=\{6\}$. Then, graphically, \\\\\\
\setlength{\unitlength}{0.5cm}
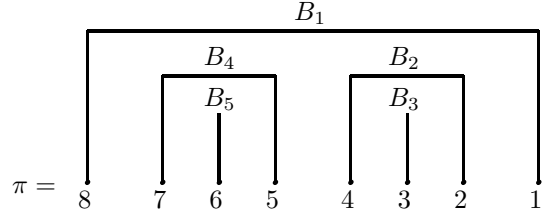
\begin{picture}(0,5)
\thicklines
\put(7,0.7){$\pi=$}
\put(9,1){\line(0,1){4}} 
\put(14.5,5.3){$B_{1}$} 
\put(12.1,4.1){$B_{4}$} 
\put(17,4.1){$B_{2}$} 
\put(12.1,3){$B_{5}$} 
\put(17,3){$B_{3}$}
 \put(8.77,0.9){\Huge.}
\put(8.5,0.3){ $8$}
 \put(11,1){\line(0,1){2.8}} 
 \put(10.77,0.9){\Huge.}
\put(10.5,0.3){ $7$}
\put(12.5,1){\line(0,1){1.8}} 
 \put(12.28,0.9){\Huge.}
\put(12,0.3){ $6$}
 \put(14,1){\line(0,1){2.8}} 
 \put(13.77,0.9){\Huge.}
\put(13.5,0.3){ $5$}
 \put(11,3.8){\line(1,0){3}} 
  \put(16,1){\line(0,1){2.8}} 
  \put(15.77,0.9){\Huge.}
\put(15.5,0.3){ $4$}
 \put(17.5,1){\line(0,1){1.8}} 
  \put(17.28,0.9){\Huge.}
\put(17,0.3){ $3$}
   \put(19,1){\line(0,1){2.8}} 
   \put(18.77,0.9){\Huge.}
\put(18.5,0.3){ $2$}
   \put(16,3.8){\line(1,0){3}} 
   \put(21,1){\line(0,1){4}} 
   \put(20.77,0.9){\Huge.}
\put(20.5,0.3){ $1$}
   \put(9,5){\line(1,0){12}} 
\end{picture}\\\\
and 
\begin{itemize}
\item $B_{1}$ is a direct predecessor of $B_{2}$ and $B_{4}$;
\item $B_2$ and $B_4$ are direct successors of $B_1$;
\item $B_{3}$ is a direct successor of $B_{2}$;
\item $B_{5}$ is a direct successor of $B_{4}$.
\end{itemize}
\end{exm}
For a given partition  $\pi\in\Pa(p)$, we define the partition $\tilde{\pi}$ as that obtained from $\pi$ by removing all singletons, $\tilde{\pi}$ will be called the \textbf{reduced partition} of $\pi$. Obviously, for $\pi\in\NC_{2}(p)$, we have $\tilde{\pi}=\pi$, but in general, for $\pi \in\NC_{2}^{1}(p)$ $\tilde{\pi}$ is a subpartition of $\pi$. For instance, \\\\if
\setlength{\unitlength}{0.5cm}
\begin{picture}(0,5)
        \thicklines
       \put(0.7,0){$\pi=$} 
      \put(2.3,0){\line(0,1){1.5}}
       \put(2.7,0){\line(0,1){1.5}}
       \put(3.1,0){\line(0,1){1}}
       \put(3.6,0){\line(0,1){1}}
        \put(2.7,1.5){\line(1,0){1.32}}
       \put(4,0){\line(0,1){1.5}}
       \put(4.5,0){\line(0,1){1.5}}
       \put(4.5,1.5){\line(1,0){1.52}}
       \put(4.8,0){\line(0,1){1}}
       \put(5.26,0){\line(0,1){0.6}}
       \put(5.7,0){\line(0,1){1}}
       \put(4.8,1){\line(1,0){0.9}}
       \put(6,0){\line(0,1){1.5}}
      \put(6.33,0){\line(0,1){1.5}}
       \put(6.8,0){$\in\NC_{2}^{1}(11)$}
      \end{picture} \hspace{6.5cm} then \hspace{3cm}\begin{picture}(0,5)
       \thicklines
       \put(-2,0){$\tilde \pi=$}
       \put(-0.5,0){\line(0,1){1.5}}
        \put(-0.5,1.5){\line(1,0){1.32}}
       \put(0.8,0){\line(0,1){1.5}}
       \put(1.3,0){\line(0,1){1.5}}
       \put(1.3,1.5){\line(1,0){1.52}}
       \put(1.6,0){\line(0,1){1}}
       \put(2.5,0){\line(0,1){1}}
       \put(1.6,1){\line(1,0){0.9}}
       \put(2.8,0){\line(0,1){1.5}}
      \put(3.3,0){$\in \NC_{2}(6).$}
      \end{picture}\\\\
\begin{df}[bm-order on noncrossing partitions with pair or singleton blocks]\label{bmorder}
Let $(\xi_{p}, \ldots, \xi_{1})$ be a given sequence of elements from a partially ordered set $(\X, \preceq)$, such that $(\xi_{p}, \ldots, \xi_{1})\sim\pi=(B_{1}, \ldots, B_{k})\in \NC_{2}^{1}(p)$, and let $\mathtt{L}: \pi\longrightarrow (\X, \preceq)$ be a label function. 
We say that the sequence $\xi:=(\xi_{p}, \ldots, \xi_{1})$ defines bm-order on the partition $\pi$ (notation $\xi \trianglelefteq \pi$) if for all $1\leq i\neq j\leq k$ the following conditions hold:
 \begin{enumerate}
  \item If $|B_{j}|=2$ and $B_{i}\prec _{\pi}B_{j}$, then $\mathtt{L}(B_{i})\preceq \mathtt{L}(B_{j})$;
  \item If $|B_{j}|=1, B_{i}\prec_{\pi} B_{j}$ and $B_{j}$ is a direct successor of $B_{i}$, then $\mathtt{L}(B_{i})=\mathtt{L}(B_{j})$.
 \end{enumerate}
We say that the sequence defines strict bm-order on the partition $\pi$ (notation $\xi\vartriangleleft \pi$) if for all $1\leq i\neq j\leq k$ the condition $(1)$ reads as: $B_{i}\prec_{\pi}B_{j}$ implies that $\mathtt{L}(B_{i})\prec \mathtt{L}(B_{j})$.
\end{df}
Note that, in a bm-ordered partition a singleton block must have the same label as its direct predecessor and pair blocks are labeled increasingly. 

Remark \ref{equalsets} suggest the following:
\begin{df}
Let $\xi:=(\xi_{p}, \ldots , \xi_{1})\in \X^p$ be a given sequence with the  adapted  partition $\pi=(B_{1}, \ldots,  B_{k})\in \NC_{2}^{1}(p)$ and $\mu:=(\mathtt{L}(B_{k}), \ldots, \mathtt{L}(B_{1}))$ be the label sequence of blocks $B_{1}, \ldots , B_{k}\in \pi$. If $\xi \trianglelefteq \pi$ (resp. $\xi \vartriangleleft \pi$), we say that the label sequence $\mu$ \textbf{defines bm-order} (resp. \textbf{defines strict bm-order}) on $\pi$ and we will use the notation $\mu \trianglelefteq _{\mathtt{L}}\pi$ (resp. $\mu \vartriangleleft_{\mathtt{L}}\pi$).
\end{df}

Our proofs of Poisson type limit theorems on discrete bm-Fock spaces eventually reduce to the study of the following combinatorial sets. 
\begin{df}
For $\pi \in \NC_{2}^{1}(p, k)$ and $\rho \in\Pi_d$ we define the following sets of sequences from $\I$:

\begin{enumerate}
\item the set 
\[
\mathrm{bmo}(\pi, \rho):=\{\xi=(\xi_{1}, \ldots, \xi_{p})\in [0, \rho]_{{\I}}^{p}: \xi\trianglelefteq \pi\}
\] 
of all sequences of elements which satisfy $0\preceq\xi_j \preceq \rho$, $\xi_j\in \I$ ($1\leq j \leq p$), and which define bm-order on $\pi$; 

\item the set 
\[
\mathrm{BMO}(\pi, \rho):=\{\mu=(\mu_{1}, \ldots, \mu_{k}): \mu_{i}=\mathtt{L}(B_{i})\in [0, \rho]_{{\I}}, \quad1\leq i\leq k \quad and \quad  \mu\trianglelefteq_{\mathtt{L}} \pi\}
\] 
of all label sequences of elements which satisfy $0\preceq\mu_j \preceq \rho$, $\mu_j\in \I$ ($1\leq j \leq k$), and define bm-order  $\pi\in \NC_{2}^{1}(p, k)$;

\item the set
\[
\mathtt{BMO}(\pi, \rho):=\{\mu=(\mu_{1}, \ldots, \mu_{k}): \mu_{i}=\mathtt{L}(B_{i})\in [0, \rho]_{{\I}}, \quad 1\leq i\leq k \quad and \quad  \mu\vartriangleleft_{\mathtt{L}} \pi\}
\]
of all label sequences of elements which satisfy $0\preceq\mu_j \preceq \rho$, $\mu_j\in \I$ ($1\leq j \leq k$), and define strict bm-order on  $\pi\in \NC_{2}^{1}(p, k)$.
\end{enumerate}

\end{df}
\begin{rem}
If $\pi\in \NC_{2}(p)$ is a pair partition, then the sets 
$\mathrm{BMO}(\pi, \rho)$ and $\mathtt{BMO}(\pi, \rho)$ are equal.
\end{rem}
We conclude this subsection by recalling the following result from \cite{bmPoisson} which plays a significant role in this study.
\begin{thm}[\cite{bmPoisson}]\label{thmp}
Let $\pi\in\NC(p, k)$ be a noncrossing partition of $p$ elements with $1\leq b(\pi)=k\leq p$ blocks, then for each positive symmetric cone $\Pi_d$ which we consider, there exists the limit 
 $$\lim_{\rho\xrightarrow[]{\Pi}\infty}\frac{|\mathtt{BMO}(\pi, \rho)|}{\vo(\rho)^{k}}=V(\pi),$$
where the function $V(\pi):=V_{\Pi_d}(\pi)$ depends on the cone $\Pi_d$ and the volume characteristic sequence $(\gamma_{n}(\Pi_d))_{n\geq 1}$, and satisfies the following recursive formula:
$$V(\pi)=\left\{
    \begin{array}{ll}
      1 & \hbox{if $ b(\pi)=1$ or  $\pi=\emptyset$ ,} \\
      \gamma_{b(\pi)}\cdot\prod\limits_{i=1}^{k}V(\pi_{i}) & \hbox{if $
    \thicklines
      \put(0,0){$\pi=$}$
      \linethickness{0.3mm}
      \put(1.5,0){\line(0,1){2.5}} 
      \put(1.6,1){ $...$}
      \put(3,0){\line(0,1){2.5}}
      \put(3.1,1){ $\pi_{1}$}
      \put(4.5,0){\line(0,1){2.5}}
      \put(4.9,1){$ ...$}
      \put(6,0){\line (0,1){2.5}}
      \put(6.1,1){ $\pi_{2}$}
      \put(7.5,0){\line (0,1){2.5}}
       \put(7.6,1){\Huge $...$}
       \put(9.2,0){\line (0,1){2.5}}
       \put(9.3,1){ $\pi_{k}$}
        \put(10.7,0){\line (0,1){2.5}}
         \put(11.1,1){$ ...$}
         \put(12.3,0){\line (0,1){2.5}}
         \put(1.5,2.5){\line(1,0){10.8}}
         \put(12.5,0){ ,}     
        
} \\
      \prod\limits_{i=1}^{k}V(\pi_{i}), & \hbox{if $\pi=\pi_1\cup \pi_2\cup\ldots \cup\pi_k$.}
    \end{array}
  \right.
$$
Here, the notation in the second case is understood as a partition $\pi$ with one outer block with $m$ element (vertical lines) and inside it there are $k$ arbitrary partitions, $\pi_{1}, \ldots, \pi_{k}$. In addition, by $\pi=\pi_1\cup \pi_2\cup\ldots \cup\pi_k$, we mean a disjoint union of sub-partitions with exactly one outer block (which could be as well a singleton).
\end{thm}

\section{Main results}

\subsection{Discrete bm-Fock space and related operators}

Let us consider a family $\{\mathcal{H}_{\xi}: \ \xi\in {\I}\subset\Pi_d \}$ of Hilbert spaces, indexed by the discrete set $\I$, with orthonormal basis $\{e_{\xi}^{m}\in \mathcal{H}_{\xi} : m=0,1, 2, 3, \cdots\}$. We assume that $\Omega:=e_{\xi}^{0}$ for all $\xi\in {\I}$ is a common unit vector, called the \emph{vacuum vector}.

Following the ideas of Muraki \cite{Mur2, Mur0}, for every $n\in {\N}$, we define
$$\varXi_{n}:=\{(\rho_{n}, \ldots, \rho_{1})\in \I^n: \quad \rho_{n}\succ \rho_{n-1}\succ \cdots \succ \rho_{1}\}.$$
The \emph{discrete bm-Fock space}, denoted by $\mathcal{F}_{bm}^{d}(\I)$, is  the Hilbert space spanned by the \emph{vacuum vector} $\Omega$ and the simple tensors of the form $h_{\rho_{n}}\otimes \cdots \otimes h_{\rho_{1}}$, where $h_{\rho}\in \mathcal{H}_{\rho}, h_{\rho}\bot\Omega$ and $(\rho_{n}, \ldots, \rho_{1})\in \varXi_{n}$; the scalar product is given by 
$$\langle h_{\rho}, h_{\eta}\rangle=0  \text{ if } \rho\neq \eta $$
and
$$\langle h_{\rho_{n}}\otimes\cdots\otimes h_{\rho_{1}}, f_{\rho_{m}}\otimes\cdots\otimes f_{\rho_{1}}\rangle=\delta_{mn}\prod_{i=1}^{n}\langle h_{\rho_{i}}, f_{\rho_{i}}\rangle,$$
where $\langle h_{\rho_{i}}, f_{\rho_{i}}\rangle$ is the scalar product in $\mathcal{H}_{\rho_{i}}$. Note that the discrete bm-Fock space $\mathcal{F}_{bm}^{d}(\I)$ is  the \emph{bm-product} \cite{JW1} of a family of Hilbert spaces $\{\mathcal{H}_{\xi}, \xi\in\I\}$, indexed by elements in  $\I\subset\Pi_d$.

In the following, we define the creation, the annihilation and the conservation operators on the discrete bm-Fock space $\mathcal{F}_{bm}^{d}(\I)$, and show their properties. Let $g_{\xi}\in\mathcal{H}_{\xi}$ be a unit vector, and let $(\rho_{n}, \ldots, \rho_{1})\in \varXi_{n}$ be an index sequence such that $h_{\rho_j}\in\mathcal{H}_{\rho_j}$ for $1\leq j \leq n$. 
\begin{enumerate} 
\item \textbf{The creation operator} $A_{g_{\xi}}^{+}$ is defined as follows:\\
  $A_{g_{\xi}}^{+}\Omega=g_{\xi}$,\\
  $A_{g_{\xi}}^{+}(h_{\rho_{n}}\otimes \cdots \otimes h_{\rho_{1}})=\left\{
                                                                      \begin{array}{ll}
                                                                        g_{\xi}\otimes h_{\rho_{n}}\otimes \cdots \otimes h_{\rho_{1}} & \hbox{if $\xi\succ\rho_{n}$,} \\
                                                                        0 & \hbox{ otherwise.}
                                                                      \end{array}
                                                                    \right.$\\                                                                   
\item \textbf{The annihilation operator} $A_{g_{\xi}}^{-}$ is defined by\\
  $A_{g_{\xi}}^{-}\Omega=0$,\\
  $A_{g_{\xi}}^{-}(h_{\rho_{n}}\otimes \cdots \otimes h_{\rho_{1}})=\left\{
                                                                      \begin{array}{ll}
                                                                        \langle g_{\xi},  h_{\rho_{n}}\rangle \cdot h_{\rho_{n-1}}\otimes \cdots \otimes h_{\rho_{1}} & \hbox{if $\xi=\rho_{n}$,} \\
                                                                        0 & \hbox{otherwise .}
                                                                      \end{array}
                                                                    \right.
  $\\
\item \textbf{The conservation operator} $A_{g_{\xi}}^{\circ}$ is defined by\\
  $A_{g_{\xi}}^{\circ}\Omega=0,$\\
  $A_{g_{\xi}}^{\circ}(h_{\rho_{n}}\otimes \cdots \otimes h_{\rho_{1}})=\left\{
                                                                          \begin{array}{ll}
                                                                            \langle g_{\xi}, h_{\rho_{n}}\rangle \cdot h_{\rho_{n}}\otimes \cdots \otimes h_{\rho_{1}} & \hbox{if $\xi=\rho_{n}$,} \\
                                                                            0 & \hbox{otherwise.}
                                                                          \end{array}
                                                                        \right.
  $\\\\
\end{enumerate}	
It follows that these operators are bounded and that $A_{g_{\xi}}^{+}$ and $A_{g_{\xi}}^{-}$ are mutually adjoint (i.e. $(A_{g_{\xi}}^{+})^{*}=A_{g_{\xi}}^{-}$). Moreover, the conservation operator $A_{g_{\xi}}^{\circ}$ is self-adjoint. In addition, we also have the following commutation relations, for the unit vectors $g_{\xi}\in\mathcal{H}_{\xi}$ and $g_{\eta}\in\mathcal{H}_{\eta}$
\begin{equation}\label{cr1}
 A_{g_{\eta}}^{+}A_{g_{\xi}}^{+}=A_{g_{\xi}}^{-}A_{g_{\eta}}^{-}=0 \quad\quad (\xi\succeq \eta)
\end{equation}
and 
\begin{equation}\label{cr3}
 A_{g_{\xi}}^{-}A_{g_{\eta}}^{\circ}=A_{g_{\xi}}^{\circ}A_{g_{\eta}}^{+}=A_{g_{\xi}}^{-}A_{g_{\eta}}^{+}=A_{g_{\xi}}^{\circ}A_{g_{\eta}}^{\circ}=0 \quad\quad ( \eta\neq \xi).
\end{equation}

For a fixed family of unit vectors $\{g_{\xi}\in \mathcal{H}_{\xi}:\ \xi \in \I \}$, we denote 
$$A_{\xi}^{\ve}:=A_{g_{\xi}}^{\ve} \text{ for } \ve\in\{\circ, -, +\}.$$

Let $\mathcal{A}$ be the $C^*$-algebra of all bounded operators on $\mathcal{F}_{bm}^{d}(\I)$ and let $\vp(A):=\langle A\Omega, \Omega \rangle$ for $A\in\mathcal{A}$ be the \emph{vacuum state}. 
For $\xi\in {\I}$, let $\mathcal{A}_{\xi}:=\text{alg}\{A_{\xi}^{+}, A_{\xi}^{-}\}$ be the $*$-algebra generated by $A_{\xi}^{+}$ and  $A_{\xi}^{-}$. In a similar way as in \cite{J.W3} one shows the following theorem. 
\begin{thm}\label{bmia}
The algebras $\{\mathcal{A}_{\xi},\ \xi\in {\I}\}$ are bm-independent in $(\mathcal{A}, \vp)$. 
\end{thm}

\subsection{Poisson type limit theorems}
In this subsection we will investigate the limits distributions of the operators $S_{\rho}(\lambda)$ defined on the discrete bm-Fock space $\mathcal{F}_{bm}^{d}(\I)$ by
\begin{equation}\label{OPD}
S_{\rho}(\lambda):=\frac{1}{\sqrt{\vo(\rho)}}\sum_{\xi\in[0, \rho]_{\I}}(A_{\xi}^{+}+A_{\xi}^{-})+\lambda\sum_{\xi\in[0, \rho]_{\I}}A_{\xi}^{\circ},
\end{equation}
as $\rho \xrightarrow[]{\Pi} \infty$ under the vacuum state $\vp$.  
To this end we define 
$$m_{p}(\lambda):=\lim_{\rho \xrightarrow []{\Pi}\infty}\vp\left((S_{\rho}(\lambda))^{p}\right),$$
the limit of $p$-th moments for such operators.\\
For each $\xi\in \Pi_d$ and $g_{\xi}\in \mathcal{H}_{\xi}$ such that $\lVert g_{\xi}\rVert=1$, as in the previous section we will use the notations $A_{\xi}^{+}:=A_{g_{\xi}}^{+}, A_{\xi}^{-}:=A_{g_{\xi}}^{-}$ and $A_{\xi}^{\circ}:=A_{g_{\xi}}^{\circ}$. Furthermore we assign the numbers $-1, 0, +1$ to the variable $\ve$ according to $\ve= -, \circ, +$.

Now we present formulation of our main results, which are the analogues of the classical Poisson limit theorem for bm-independent random variables, indexed by elements of positive symmetric cones. 
\begin{thm}\label{mthm}
 For any positive integer $p\geq 2$, one has
 \begin{equation}\label{mp}
 m_{p}(\lambda)=\sum_{\pi\in \NC_{2}^{1, i}(p)}\lambda^{s(\pi)}\cdot  V(\tilde \pi),
 \end{equation}
 where 
\begin{equation}\label{bmolim}
V(\tilde \pi)=\lim\limits_{\rho\xrightarrow[]{\Pi} \infty}\frac{|\mathrm{BMO}(\pi, \rho)|}{(\vo(\rho))^{b(\tilde \pi)}}.
\end{equation} 
 Here  $b(\tilde \pi)$ is the number of blocks in the reduced partition $\tilde \pi$ of $\pi $ and the limit \eqref{bmolim} exists for every partition $\pi\in \NC_{2}^{1, i}(p)$. The combinatorial function $V(\tilde  \pi):=V_{\Pi_d}(\tilde{\pi})$ depends on the positive cone $\Pi_d$ and its volume characteristic $\gamma_{n}(\Pi_d)$. Moreover, it is multiplicative and satisfies the following  recursive formula
 \begin{equation}\label{rfv}
 V(\tilde \pi)=\left\{
    \begin{array}{ll}
      1 & \hbox{if $b(\tilde \pi)=1$ or  $\tilde \pi=\emptyset$} \\
      \gamma_{|\tilde \pi |}(\Pi_d)\cdot V(\pi') & \hbox{if $
     
      \put(0.2,0){$\tilde \pi=$}$
      \begin{picture}(5,2)
       \thicklines
      \put(1.3,0){\line(0,1){1.5}} 
      \put(3,0.2){$\pi'$}
      \put(5,0){\line(0,1){1.5}}
      \put(1.3,1.5){\line(1,0){3.72}}
      \end{picture}} \\
      \prod\limits_{i=1}^{k}V(\pi_{i}) & \hbox{if \hspace{0.14cm}$\tilde \pi=\pi_{1}\cup \pi_{2}\cup \cdots\cup\pi_{k}$,}
    \end{array}
  \right.
\end{equation} 
\end{thm}

\subsection{Examples and remarks}

The following examples and figures explain the rules for the calculations of $V(\tilde \pi)$ for a given partition $\pi\in\NC_{2}^{1, i}(p)$. For example, consider the following  partition $\pi$:
 
$$\pi=\{\{1, 15\}, \{2\}, \{3, 9\}, \{4, 8\}, \{5\}, \{6\}, \{7\}, \{10, 13\}, \{11\},  \{12\}, \{14\}\}\in\NC_{2}^{1, i}(15)$$ Then we can draw the following pictures of $\pi$ and $\tilde{\pi}$\\\\
      \begin{picture}(4.8,2)
       \thicklines
       \put(-0.1,0){$\pi=$}
      \put(1.3,0){\line(0,1){2}} 
      \put(1.5,0){\line(0,1){1.5}}
       \put(1.7,0){\line(0,1){1.5}}
       \put(2.1,0){\line(0,1){1}}
       \put(2.6,0){\line(0,1){1}}
        \put(1.7,1.5){\line(1,0){1.32}}
       \put(3,0){\line(0,1){1.5}}
       \put(3.5,0){\line(0,1){1.5}}
       \put(3.5,1.5){\line(1,0){1.52}}
       \put(3.8,0){\line(0,1){1}}
       \put(4,0){\line(0,1){0.6}}
       \put(4.26,0){\line(0,1){0.6}}
       \put(4.5,0){\line(0,1){0.6}}
       \put(4.7,0){\line(0,1){1}}
       \put(3.8,1){\line(1,0){0.9}}
       \put(5,0){\line(0,1){1.5}}
      \put(5.5,0){\line(0,1){2}}
      \put(5.25,0){\line(0,1){1.5}}
      \put(1.3,2){\line(1,0){4.22}}
      \end{picture} \hspace{2.5cm} and \hspace{2cm}\begin{picture}(4.8,2)
       \thicklines
       \put(-0.1,0){$\tilde \pi=$}
      \put(1.3,0){\line(0,1){2}} 
       \put(1.7,0){\line(0,1){1.5}}
        \put(1.7,1.5){\line(1,0){1.32}}
       \put(3,0){\line(0,1){1.5}}
       \put(3.5,0){\line(0,1){1.5}}
       \put(3.5,1.5){\line(1,0){1.52}}
       \put(3.8,0){\line(0,1){1}}
       \put(4.7,0){\line(0,1){1}}
       \put(3.8,1){\line(1,0){0.9}}
       \put(5,0){\line(0,1){1.5}}
      \put(5.5,0){\line(0,1){2}}
      \put(1.3,2){\line(1,0){4.22}}
      \end{picture}\\\\
Recall that, for $\Pi_2={\R}_{+}^{2}$ and $\Pi_1=\Lambda_{1}^{1}$, the volume characteristic sequence is $\gamma_{n}=\frac{1}{n^{2}}$, and for $\Pi_3={\R}_{+}^{3}$, $\gamma_{n}=\frac{1}{n^{3}}$, while for $\Pi_2=\Lambda_{2}^{1}$ and $\Pi_2=\Symm$, $\gamma_{n}=\frac{24}{(3n-1)3n(3n+1)}$. Now we can compute $V(\tilde \pi)$ for such cones using the recursive formula \eqref{rfv}:\\\\
$1.\quad\Pi_2={\R}_{+}^{2}$ and $\Pi_1=\Lambda_{1}^{1}.$ \\\\
 V($\tilde \pi$)=V(\hspace{-0.5cm}\begin{picture}(4.8,2)
       \thicklines
      \put(1.3,0){\line(0,1){2}} 
       \put(1.7,0){\line(0,1){1.5}}
        \put(1.7,1.5){\line(1,0){1.32}}
       \put(3,0){\line(0,1){1.5}}
       \put(3.5,0){\line(0,1){1.5}}
       \put(3.5,1.5){\line(1,0){1.52}}
       \put(3.8,0){\line(0,1){1}}
       \put(4.7,0){\line(0,1){1}}
       \put(3.8,1){\line(1,0){0.9}}
       \put(5,0){\line(0,1){1.5}}
      \put(5.5,0){\line(0,1){2}}
      \put(1.3,2){\line(1,0){4.22}}
      \end{picture}\hspace{0.5cm})=$\frac{1}{16}\cdot$ V\hspace{0.1cm}(\begin{picture}(5,2)
       \thicklines
    
       \put(0.19,0){\line(0,1){1.5}}
        \put(0.19,1.5){\line(1,0){1.32}}
       \put(1.49,0){\line(0,1){1.5}}
    
      \end{picture}\hspace{-1.65cm})$\cdot V(\begin{picture}(5,2)
       \thicklines
       \put(0.2,0){\line(0,1){1.5}}
       \put(0.2,1.5){\line(1,0){1.52}}
       \put(0.5,0){\line(0,1){1}}
       \put(1.4,0){\line(0,1){1}}
       \put(0.5,1){\line(1,0){0.9}}
       \put(1.7,0){\line(0,1){1.5}}
      \end{picture}$\hspace{-1.55cm})=$\frac{1}{16}\cdot 1\cdot \frac{1}{4}\cdot V(\begin{picture}(5,2)
       \thicklines
       \put(0.2,0){\line(0,1){1}}
       \put(1.1,0){\line(0,1){1}}
       \put(0.2,1){\line(1,0){0.9}}
      \end{picture}\hspace{-1.85cm})=\frac{1}{64}$
\\\\ $2.\quad\Pi_3={\R}_{+}^{3}.$\\\\
V($\tilde \pi$)=V(\hspace{-0.5cm}\begin{picture}(4.8,2)
       \thicklines
      \put(1.3,0){\line(0,1){2}} 
       \put(1.7,0){\line(0,1){1.5}}
      
        \put(1.7,1.5){\line(1,0){1.32}}
       \put(3,0){\line(0,1){1.5}}
       \put(3.5,0){\line(0,1){1.5}}
       \put(3.5,1.5){\line(1,0){1.52}}
       \put(3.8,0){\line(0,1){1}}
       \put(4.7,0){\line(0,1){1}}
       \put(3.8,1){\line(1,0){0.9}}
       \put(5,0){\line(0,1){1.5}}
      
      \put(5.5,0){\line(0,1){2}}
      \put(1.3,2){\line(1,0){4.22}}
      \end{picture}\hspace{0.5cm})=$\frac{1}{64}\cdot$ V(\begin{picture}(5,2)
       \thicklines
    
       \put(0.2,0){\line(0,1){1.5}}
        \put(0.2,1.5){\line(1,0){1.32}}
       \put(1.5,0){\line(0,1){1.5}}
    
      \end{picture}\hspace{-1.65cm})$\cdot V(\begin{picture}(5,2)
       \thicklines
       \put(0.2,0){\line(0,1){1.5}}
       \put(0.2,1.5){\line(1,0){1.52}}
       \put(0.5,0){\line(0,1){1}}
       \put(1.4,0){\line(0,1){1}}
       \put(0.5,1){\line(1,0){0.9}}
       \put(1.7,0){\line(0,1){1.5}}
      \end{picture}$\hspace{-1.55cm})=$\frac{1}{64}\cdot 1\cdot \frac{1}{8}\cdot V(\begin{picture}(5,2)
       \thicklines
       \put(0.2,0){\line(0,1){1}}
       \put(1.1,0){\line(0,1){1}}
       \put(0.2,1){\line(1,0){0.9}}
      \end{picture}\hspace{-1.85cm})=\frac{1}{512}$\\\\
$3.\quad\Pi_2=\Lambda_{2}^{1}$ and $\Pi_2=\Symm$.\\\\
V($\tilde \pi$)=V(\hspace{-0.5cm}\begin{picture}(5,2)
       \thicklines
      \put(1.3,0){\line(0,1){2}} 
       \put(1.7,0){\line(0,1){1.5}}
        \put(1.7,1.5){\line(1,0){1.32}}
       \put(3,0){\line(0,1){1.5}}
       \put(3.5,0){\line(0,1){1.5}}
       \put(3.5,1.5){\line(1,0){1.52}}
       \put(3.8,0){\line(0,1){1}}
       \put(4.7,0){\line(0,1){1}}
       \put(3.8,1){\line(1,0){0.9}}
       \put(5,0){\line(0,1){1.5}}
      \put(5.5,0){\line(0,1){2}}
      \put(1.3,2){\line(1,0){4.22}}
      \end{picture}\hspace{0.4cm})=$\frac{2}{143}\cdot$ V(\begin{picture}(5,2)
       \thicklines
    
       \put(0.2,0){\line(0,1){1.5}}
        \put(0.2,1.5){\line(1,0){1.32}}
       \put(1.5,0){\line(0,1){1.5}}
    
      \end{picture}\hspace{-1.65cm})$\cdot V(\begin{picture}(5,2)
       \thicklines
       \put(0.2,0){\line(0,1){1.5}}
       \put(0.2,1.5){\line(1,0){1.52}}
       \put(0.5,0){\line(0,1){1}}
       \put(1.4,0){\line(0,1){1}}
       \put(0.5,1){\line(1,0){0.9}}
       \put(1.7,0){\line(0,1){1.5}}
      \end{picture}$\hspace{-1.55cm})=$\frac{2}{143}\cdot 1\cdot \frac{12}{105}\cdot V(\begin{picture}(5,2)
       \thicklines
       \put(0.2,0){\line(0,1){1}}
       \put(1.1,0){\line(0,1){1}}
       \put(0.2,1){\line(1,0){0.9}}
      \end{picture}\hspace{-1.85cm})=\frac{8}{505}$\\\\
For illustration, the first few moments ($m_{1}(\lambda), \ldots, m_{6}(\lambda)$) of the monotone Poisson distribution (cf.\ \cite{Mur0}) in comparison with some of the bm-Poisson analogues are listed below (observe that always $m_{1}(\lambda)=0, m_{2}(\lambda)=1$ and $m_{3}(\lambda)=\lambda$). 
\vspace{1cm}
\begin{center}
\hspace{0.5cm}
\begin{tabular}{ c c c c c c c} 

\hspace{-2cm}Monotone Poisson distribution &&&& \hspace{2cm}bm-case for $\Pi_3={\R}_{+}^{3}$ \\
\hspace{-2cm}$0$ \vspace{0.2cm}&&&&\hspace{2cm} $0$\\
\hspace{-2cm}$1$ \vspace{0.2cm}&&&&\hspace{2cm} $1$\\
\hspace{-2cm}$\lambda$ \vspace{0.2cm}&&&& \hspace{2cm}$\lambda$ \\ 
\hspace{-2cm}$\lambda^{2}+\frac{3}{2}$ \vspace{0.2cm}&&&&\hspace{2cm}$\lambda^{2}+\frac{9}{8}$\\
\hspace{-2cm}$\lambda^{3}+\frac{7}{2}\lambda$ \vspace{0.2cm}&&&& \hspace{2cm}$\lambda^{3}+\frac{19}{8}\lambda$ \\ 
\hspace{-2cm}$\lambda^{4}+\frac{9}{2}\lambda^{2}+\frac{5}{2}$  \vspace{2cm} &&&& \hspace{2cm}$\lambda^{4}+\frac{15}{4}\lambda^{2}+\frac{31}{24}$ \\ 
\end{tabular}
\end{center}

\vspace{0.2cm}
\begin{center}
\hspace{0.5cm}
\begin{tabular}{ c c c c c c c} 
\hspace{-1cm} bm-cases for $\Pi_2=\R_+^2$ and $\Pi_1=\Lambda_{1}^{1}$&&&& \hspace{0.2cm}bm-cases for $\Pi_2=\Lambda_{2}^{1}$ and $\Pi_2=\Symm$ \\

\hspace{-1cm}$0$ \vspace{0.2cm}&&&&\hspace{0.2cm} $0$\\
\hspace{-1cm}$1$ \vspace{0.2cm}&&&& \hspace{0.2cm} $1$ \\ 
\hspace{-1cm}$\lambda$ \vspace{0.2cm}&&&&\hspace{0.2cm}$\lambda$\\
\hspace{-1cm}$\lambda^{2}+\frac{5}{4}$ \vspace{0.2cm}&&&& \hspace{0.2cm}$\lambda^{2}+\frac{39}{35}$ \\ 
\hspace{-1cm}$\lambda^{3}+\frac{11}{4}\lambda$ \vspace{0.2cm}&&&& \hspace{0.2cm}$\lambda^{3}+\frac{12}{35}\lambda$ \\ 
\hspace{-1cm}$\lambda^{4}+\frac{9}{2}\lambda^{2}+\frac{59}{36}$ \vspace{0.2cm} &&&& \hspace{0.2cm}$\lambda^{4}+\frac{129}{35}\lambda^{2}+\frac{443}{350}$ \\

\end{tabular}
\end{center}
\vspace{1cm}
In the special case $\lambda=0$, the operators \eqref{OPD} become  
$$S_{\rho}(0)=\frac{1}{\sqrt{\vo(\rho)}}\sum_{\xi\preceq \rho, \hspace{0.1cm}\xi\in\I}(A_{\xi}^{+}+A_{\xi}^{-}).$$
Since $\{A_{\xi}^{+}+A_{\xi}^{-}, \xi\in{\I}\}$ are bm-independent on $(\mathcal{A}, \vp)$, then we obtain a new example of noncommutative ``de Moivre-Laplace theorem" which is a kind of central limit theorem. Namely the limit distribution of the operator $S_{\rho}(0)$ at $\rho\xrightarrow[]{\Pi} \infty$ under the vacuum state $\vp$ is given by 
$$m_{2n+1}(0)=0\quad\text{ and }\quad m_{2n}(0)=\sum_{\pi\in\NC_{2}(2n)}V(\pi).$$
On the other hand, by the bm-central limit theorem \cite{J.W3} the  even moments  $m_{2n}(0)$ can be written as follows 
$$m_{2n}(0)=g_{n}=\sum_{k=1}^{n}\gamma_{k}(\Pi_d)g_{k-1}g_{n-k}, \quad\quad g_{0}=g_{1}=1.$$
Hence, we have another combinatorial description for the bm-central limit theorem associated with symmetric cones \cite{J.W3}.\\
It is worthwhile to mention that it is also possible to extend the results of this paper to two classes of positive cones studied in \cite{OW21019} which are non-symmetric. Namely, the sectorial cones 
$$\Omega_{\bm u}^{n}:=\biggl\{\sum_{j=1}^{n}a_{j}u_{j}: \bm u:=(u_{1}, \ldots, u_{n}), a_{j}\geq 0, u_{j}\in\mathbb{R}^{n} \text{ and } u_{1}, \ldots, u_{n} \text{ are linearly independent} \biggr\}$$
and the circular cones
$$C_{\theta}^{n}:=\{(t; x)\in \mathbb{R}_{+}\times \mathbb{R}^{n}: ||x||\leq t\cdot \tan \theta\}.$$

In what follows we will consider the subspace of $\mathcal{F}_{bm}^{d}(\I)$
$$\mathcal{F}_{0}:=\text{ span }\{\Omega, g_{\xi_{n}}\otimes \cdots \otimes g_{\xi_{1}}: \xi_{n}\succ \xi_{n-1}\succ \cdots \succ \xi_{1},\ n \in \N\},
$$
and by $\B_{\xi}^{\ve}:=A_{\xi}^{\ve}|_{_{\mathcal{F}_{0}}}$ for $\ve\in\{-, \circ, +\}$ and $\xi\in \I$ we will denote the restrictions of $A_{\xi}^{\ve}$ to the subspace $\mathcal{F}_{0}$. Observe that each operator $A_{\xi}^{\ve}$  preserves the subspace $\mathcal{F}_{0}$, and moreover 
\begin{equation}\label{eqsub}
\B_{\xi}^{\circ}=\B_{\xi}^{+}\B_{\xi}^{-} \text { is in } \mathcal{F}_{0}.
\end{equation}
Let $\mathcal{C}$ be the unital $*$-algebra of all bounded operators on $\mathcal{F}_{0}$. Hence, if we define $\mathcal{C}_{\xi}$ to be the $*$-algebra generated by $\B_{\xi}^{+}, \B_{\xi}^{-}$ and  $\B_{\xi}^{\circ}$, we get the following crucial Corollary.
\begin{cor}\label{bmc}
The algebras $\{\mathcal{C}_{\xi}, \xi\in {\I}\}$ are bm-independent with respect to $\vp$ in $(\mathcal{C}, \vp).$
\end{cor}
\subsection{Proof of Theorem \ref{mthm}}
The proof of Theorem \ref{mthm} consists of several reductions and combinatorial considerations. First we will show that the limit \eqref{mp} can be reduced combinatorially to a sum over noncrossing partitions $\pi$ with pair or inner singleton blocks only. Next, we shall prove that the only terms which survive in the limit are those which satisfy the bm-ordered labellings $\mu\trianglelefteq_{\mathtt{L}}\pi$ for a partition $\pi$. Finally, we shall show that the limit when $\rho \xrightarrow[]{\Pi}\infty$ is the ratio of the cardinality $|\mathrm{BMO}(\pi, \rho)|$ and the Euclidean volume $\vo(\rho)^{b(\tilde{\pi})}$ for $\pi\in\NC_{2}^{1, i}(p)$, where $\tilde{\pi}$ is the reduced partition of $\pi$. The rest of the proof then goes in a similar manner as for Theorem 6.1 in \cite{bmPoisson}.
\subsubsection{Combinatorial reduction}
We start the proof observing that $\vp((S_{\rho}(\lambda))^{p})$ can be written as 
\begin{align}\label{credu1}
 \vp((S_{\rho}(\lambda))^{p})&=\frac{1}{\vo(\rho)^{\frac{p}{2}}}\sum_{(\ve_{p}, \ldots, \ve_{1})\in \{-1, 0, +1\}^{p}}\hspace{0.2cm}\sum_{(\xi_{p}, \ldots, \xi_{1})\in [0, \rho]_{\I}^{p}}\vp(A_{\xi_{p}}^{\ve_{p}} \cdots A_{\xi_{1}}^{\ve_{1}})\lambda_{\rho}^{\sum_{i=1}^{p}\delta_{0}(\ve_{i})}\nonumber\\
 &=\frac{1}{\vo(\rho)^{\frac{p}{2}}}\sum_{(\ve_{p}, \ldots, \ve_{1})\in \{-1, 0, +1\}^{p}}\hspace{0.2cm}\sum_{(\xi_{p}, \ldots, \xi_{1})\in [0, \rho]_{\I}^{p}}\vp(\B_{\xi_{p}}^{\ve_{p}} \cdots \B_{\xi_{1}}^{\ve_{1}})\lambda_{\rho}^{\sum_{i=1}^{p}\delta_{0}(\ve_{i})}
\end{align} 
where $$\delta_{0}(\ve_{i})=\left\{
                                                                                                \begin{array}{ll}
                                                                                                  1 & \hbox{if $\ve_{i}=0$} \\
                                                                                                  0 & \hbox{otherwise,}
                                                                                                \end{array}
                                                                                              \right.$$
and $\lambda_{\rho}=\lambda\cdot \sqrt{\vo(\rho)}.$\\
For $p\in {\N}$, we will use the following notations $\bm{\ve}:=(\ve_{p},\ldots, \ve_{1})\in \{-1, 0, +1\}^{p}$, $\bm{\xi}:=(\xi_{p}, \ldots, \xi_{1})\in [0, \rho]_{{\I}}^{p}$ and $\B_{\bm\xi}^{\bm \ve}:=\B_{\xi_{p}}^{\ve_{p}} \cdots \B_{\xi_{1}}^{\ve_{1}}$. By a simple induction argument on $p$, one obtains the following lemma which describes a necessary conditions for non-vanishing of the mixed moments $\vp(\B_{\bm \xi}^{\bm \ve})$. The reader is refereed to \cite{CGW2021} for similar arguments in the setting of weakly monotone Fock space and to \cite{CVLYG2007} in the setting of interacting Fock space.
\begin{lemm}\label{lem1}
 For a positive integer $p\geq 2$, let $\bm{\ve}\in \{-1, 0, +1\}^{p}$ and $\bm \xi\in [0, \rho]_{{\I}}^{p}$.\\
 If $\vp(\B_{\bm \xi}^{\bm \ve})\neq 0$, then the sequence $\bm \ve$ satisfies the following three conditions:
 \begin{enumerate}
 \item $\ve_{1}=+1, \ve_{p}=-1,$
  \item $\sum\limits_{i=1}^{p}\ve_{i}=0,$
  \item $\sum\limits_{i=1}^{k}\ve_{i}\geq 0$\quad  for $ k=1, \ldots, p-1.$
 \end{enumerate}
\end{lemm}
From now on, by ${\D}_{p}$ we denote the set of all sequences $\bm \ve=(\ve_{p}, \ldots, \ve_{1})\in \{-1, 0, +1\}^{p}$ satisfying the conditions (1)-(3) in the statement of Lemma \ref{lem1}.
According to the same Lemma and using the above  notations,~\eqref{credu1} becomes 
\begin{equation}\label{credu2}
 \vp((S_{\rho}(\lambda))^{p})=\frac{1}{\vo(\rho)^{\frac{p}{2}}}\sum_{\bm \ve\in \D_{p}}\hspace{0.2cm}\sum_{\bm \xi\in [0, \rho]_{{\I}}^{p}}\vp(\B_{\bm\xi}^{\bm\ve})\lambda_{\rho}^{\sum_{i=1}^{p}\delta_{0}(\ve_{i})}.
\end{equation}
For a fixed sequence $\bm{\ve}\in \D_{p}$ and for each $1\leq k\leq p-1$ ($p\geq 2$) such that $\ve_{k}=+1$, by (2) and (3) in Lemma \ref{lem1} one can define
\begin{equation}\label{tk}
 T(k):=\min \{l: l>k \quad s.t.\quad  \ve_{k}+\cdots+\ve_{l}=0\}.
\end{equation}
By definition, $\ve_{T(k)}=-1$. Furthermore, if there exists $k<m<T(k)$ such that $\ve_{m}=+1$, then 
$$\ve_{k}+\cdots+\ve_{m-1}>0 \text{ and } \ve_{m}+\ve_{m+1}+\cdots+\ve_{T(k)}<0.$$
This implies that  
$$k<m<T(m)<T(k).$$
Therefore, if we define a partition 
$$\pi:=\{\{k, T(k)\}: \ve_{k}=+1\} \cup\{\{j\}: \ve_{j}=0\},$$
we obtain $\pi\in \NC_{2}^{1}(p)$. In other words, for any sequence $\bm \ve\in \D_{p}$ we can uniquely associate a noncrossing partition with blocks being a pair or a singleton, and we will often use the identification $\D_{p}\ni \bm\ve\equiv \pi\in\NC_{2}^{1}(p)$.
\begin{lemm}\label{lem2}
For $p\geq 2$, let $\bm \ve\in \D_{p}, \bm\xi\in [0, \rho]_{{\I}}^{p}$ and $1\leq k\leq p-1$ such that $\ve_{k}=+1$.
 If $\vp(\B_{\bm{\xi}}^{\bm{\ve}})\neq 0$ then $\xi_{k}=\xi_{T(k)}.$ 
\end{lemm}
\bproof
We use induction over $p\geq 2$:\\
For $p=2$, we have just one pair block $\{1, T(1)=2\}$, then \eqref{cr3} gives the thesis.\\
For $p=3$, we have the sequences $\bm\ve=(-1, 0, +1)$ and $\bm\xi=(\xi_{3},\xi_{2},\xi_{1})$. Hence, $T(1)=3$ and the assumption $\vp(\B_{\xi}^{\ve})\neq 0$ yields $\xi_{1}=\xi_{2}=\xi_{3}=\xi_{T(1)}$.\\
We assume that the lemma holds for all $r<p$ and that $\vp(\B_{\bm \xi}^{\bm \ve})\neq 0$. Consider 
$$\text{$1\leq k_{0}\leq p-1$ such that $ \ve_{k_{0}}=+1$  and $T(k_{0})\leq T(k)$ for all $ 1\leq k\leq p-1$ with $\ve_{k}=+1$},$$
that is, $$T(k_{0}):=\min \{ T(k): 1\leq k\leq p-1\}.$$\\
Then, if $|T(k_{0})-k_{0}|\geq 2$, we have $\ve_{j}=0$ for all $k_{0}<j<T(k_{0})$, and if there exists  $1\leq i\leq k_{0}$, then $\ve_{i}\in\{0, +1\}$. Suppose initially that $T(k_{0})=p$, hence $k_{0}=1$. In this case,
$$\vp(\B_{\bm\xi}^{\bm\ve})=\vp(\B_{\xi_{p}}^{-}\B_{\xi_{p-1}}^{\circ}\cdots \B_{\xi_{2}}^{\circ}\B_{\xi_{1}}^{+})\neq 0.$$
Hence, 
$$\xi_{1}=\xi_{2}=\cdots=\xi_{p-1}=\xi_{p}=\xi_{T(1)}.$$
If instead $T(k_{0})< p$, the assumption $\vp(\B_{\bm\xi}^{\bm\ve})\neq 0$ yields
\begin{equation*}
 \B_{\xi_{k_{0}}}^{+}(\B_{\xi_{k_{0}-1}}^{\ve_{k_{0}-1}}\cdots \B_{\xi_{1}}^{+})\Omega\neq 0 ,
\end{equation*}
thus
\begin{equation*}
\B_{\xi_{k_{0}}}^{+}(\B_{\xi_{k_{0}-1}}^{\ve_{k_{0}-1}}\cdots \B_{\xi_{1}}^{+})\Omega=\B_{\xi_{k_{0}}}^{+}w=g_{\xi_{k_{0}}}\otimes w\neq 0,
\end{equation*}
where 
$$w=g_{\xi_{s_{l}}}\otimes \cdots \otimes g_{\xi_{s_{1}}}\otimes g_{\xi_{1}} \text{ and 
} \xi_{s_{1}}, \ldots, \xi_{s_{l}}\in \{\xi_{2}, \ldots, \xi_{k_{0}-1}\} \text{ with }\xi_{k_{0}}\succ \xi_{s_{l}}\succ \cdots \succ \xi_{s_{1}}\succ \xi_{1}.$$
Consequently, since
$$\B_{\xi_{T(k_{0})-1}}^{\circ}\cdots \B_{\xi_{k_{0}+1}}^{\circ}g_{\xi_{k_{0}}}\otimes w\neq 0\Longrightarrow \xi_{T(k_{0})-1}=\cdots =\xi_{k_{0}+1}=\xi_{k_{0}},$$
and 
$$\B_{\xi_{T(k_{0})}}^{-}(\B_{\xi_{T(k_{0})-1}}^{\circ}\cdots \B_{\xi_{k_{0}+1}}^{\circ})g_{\xi_{k_{0}}}\otimes w=\B_{\xi_{T(k_{0})}}^{-}(\B_{\xi_{T(k_{0})}}^{\circ})^{n}g_{\xi_{k_{0}}}\otimes w=\B_{\xi_{T(k_{0})}}^{-}g_{\xi_{k_0}}\otimes w\neq 0,$$
where $n=T(k_{0})-1-k_{0}$.\\
Therefore, $\xi_{T(k_{0})}=\xi_{k_{0}} \text{ and } \vp(\B_{\bm \xi}^{\bm \ve})=\vp(\B_{\bm{\xi}'}^{\bm{\ve}'})$, 
where
$$\bm {\ve}'=(\ve_{p}, \ldots, \ve_{T(k_{0})+1}, \ve_{k_{0}-1}, \ldots, \ve_{1}) \text{ and } \bm{\xi}'=(\xi_{p}, \ldots, \xi_{T(k_{0})+1}, \xi_{k_{0}-1}, \ldots, \xi_{1}).$$
Since $r=p-(n+2)<p$, the thesis follows after exploiting the induction assumption for the sequences $\bm{\ve}'$ and $ \bm{\xi}'$.
\eproof
\begin{cor}\label{cor1}
Let $\bm{\ve}\in\D_{p}$ and $\bm{\xi}\in[0, \rho]_{\I}^{p}$. If $\vp(\B_{\bm\ve}^{\bm\xi})\neq 0$, then the singletons in the partition $\pi=\{\{k, T(k)\}, \ve_{k}=+1\}\cup\{\{j\}, \ve_{j}=0\}$ which is associated with the sequence $\bm \ve$, cannot be outer blocks, i.e. $\pi\in \NC_{2}^{1, i}(p)$.
\end{cor}
Using Lemma \ref{lem2} and Corollary \ref{cor1}, the identity \eqref{credu2} becomes:
\begin{equation}\label{credu3}
  \vp((S_{\rho}(\lambda))^{p})=\frac{1}{\vo(\rho)^{\frac{p}{2}}}\sum_{\pi\in \NC_{2}^{1, i}(p)}\hspace{0.1cm}\sum_{[0, \rho]_{{\I}}^{p}\ni \bm \xi\sim\pi}\vp(\B_{\bm \xi}^{\bm \ve})\lambda^{s(\pi)},
\end{equation}
 where the second summation is over all sequence $\bm \xi=(\xi_{p}, \ldots, \xi_{1})\in[0, \rho]_{{\I}}^{p}$ associated with a given partition $\pi\in \NC_{2}^{1, i}(p)$.
\subsubsection{Reduction to bm-ordered noncrossing partitions with pair or singleton blocks:}
In the next lemma we will show  that the only terms which survive in the limit \eqref{mp} come from those sequences establishing a bm-order on noncrossing partitions with pair or inner singleton blocks, i.e., those satisfying the two conditions in Definition \ref{bmorder}.
\begin{lemm}[bm-order] 
Let $\pi\in \NC_{2}^{1, i}(p)$ be a noncrossing partition with pair or inner singleton blocks given by the sequences $\bm \ve\in \D_{p}$ and  $\bm{\xi}\in[0, \rho]_{\I}^{p}$ such that $\bm{\xi}\sim\pi$. If $\vp(\B_{\bm \xi}^{\bm \ve})\neq 0$ then the sequence $\bm \xi$ establishes a bm-order on $\pi$, i.e. $\bm \xi \unlhd \pi $.
\end{lemm}
\bproof
By induction on $p\geq 2$:\\
For $p=2$, we have $\bm \ve=(-1, +1), \bm\xi=(\xi_{2}, \xi_{1})$ and  $\pi=\{\{1, 2\}\}$. Then $\vp(\B_{\bm \xi}^{\bm \ve})\neq 0$ implies that  $\xi_{1}=\xi_{2}$ and $(\xi_{2}, \xi_{1})$ establishes bm-order on $\pi=\{\{1, 2\}\}$.\\
For $p=3$, our sequences are $\bm\ve=(-1, 0, +1)$ and $\bm \xi=(\xi_{3}, \xi_{2}, \xi_{1})$ with the associated partition $\pi=\{\{1, 3\}, \{2\}\}$. Then
$$\vp(\B_{\bm\xi}^{\bm\ve})\neq 0\Longrightarrow \B_{\xi_{3}}^{-}\B_{\xi_{2}}^{0}\B_{\xi_{1}}^{+}\Omega\neq 0\Longrightarrow \xi_{1}=\xi_{2}=\xi_{3},$$
and $\bm\xi=(\xi_{3}, \xi_{2}, \xi_{1})$ establishes bm-order on $\pi=\{\{1, 3\},\{2\}\}$.\\
Let us consider $\bm\xi=(\xi_{p}, \ldots, \xi_{1})\sim \pi\in \NC_{2}^{1, i}(p)$ and $\bm{\ve}=(\ve_{p}, \ldots, \ve_{1})\in \D_{p}$. We assume that the statement is true for each $q<p$ and $\vp(\B_{\bm\xi}^{\bm\ve})\neq 0$. We have the following two cases, depending on the position of $T(1)$. 
\begin{enumerate}
\item $T(1)<p:$\\
Denote $C:=\B_{\xi_{T(1)}}^{\ve_{T(1)}} \cdots \B_{\xi_{1}}^{+}$ and $D:=\B_{\xi_{p}}^{-} \cdots \B_{\xi_{T(1)+1}}^{\ve_{T(1)+1}}$. As $\xi_{1}=\xi_{T(1)}$ and $\ve_{T(1)}=-1$, one gets 
$$C\Omega=\Omega \text{ and } 0\neq \vp(\B_{\bm \xi}^{\bm \ve})=\vp(D).$$
Then $\ve_{T(1)+1}=+1$.\\
If we denote by $\bm{\ve}':=(\ve_{T(1)}, \ldots, \ve_{1}), \bm{\xi}':=(\xi_{T(1)}, \ldots, \xi_{1}), \bm {\ve}'':=(\ve_{p}, \ldots, \ve_{T(1)+1})$ and $\bm{\xi}'':=(\xi_{p}, \ldots, \xi_{T(1)+1})$, then  $\bm{\ve}'\equiv\pi'\in \NC_{2}^{1, i}(T(1))$ and $\bm{\ve}''\equiv \pi''\in \NC_{2}^{1, i}(p-T(1))$, where $\bm{\xi}'\sim\pi'$ and $\bm{\xi}''\sim\pi''$.\\
Hence, 
$$1=\vp(C)=\vp(\B_{\bm \xi'}^{\bm\ve'})\neq 0 \quad and \quad \vp(D)=\vp(\B_{\bm \xi''}^{\bm \ve''})\neq 0.$$
Since $T(1)<p$ and $p-T(1)<p$, we can use the induction hypothesis for $\vp(\B_{\bm\xi'}^{\bm\ve'})\neq 0$ to obtain that $\bm\xi'$ establishes bm-order on $\pi'$, and for $\vp(\B_{\bm\xi''}^{\bm\ve''})\neq 0$ to obtain that $\bm \xi''$ establishes bm-order on $\pi''$.\\
However, there is no relation between the sequence $\bm \xi'\sim\pi'$ and the sequence $\bm \xi''\sim\pi''$, i.e. for $\eta\in \bm\xi'$ and $\mu\in \bm \xi''$ all three possibilities $\eta\succeq\mu,\eta\preceq \mu$ and $\eta \nsim \mu$ are allowed. Therefore, $\bm\xi$ establishes bm-order on $\pi=\pi'\cup \pi''$.
\item $T(1)=p:$\\
As $\ve_{p}=-1$, let us take $j=\min\{i: 2\leq i \leq p, \ve_{i}=-1\}$. Note that if $j=2$, then $T(1)=2=p$, and if $j=p$ (we have one outer pair block $\{1, p\}$ and $(p-2)$ inner singletons), then $\ve_{i}=0$ and $ \xi_{i}=\xi_{1}=\xi_{T(1)}$ for $i=2, \ldots, p-1$. In both cases, the sequence $\bm\xi$ establishes bm-order on $\pi\in \NC_{2}^{1, i}(p)$, and 
$$\vp(\B_{\bm\xi}^{\bm\ve})=\vp(\B_{\xi_{p}}^{-}\B_{\xi_{p-1}}^{\circ}\cdots B^{\circ}_{\xi_{\xi_{2}}}\B_{\xi_{1}}^{+})=\vp(\B_{\xi_{p}}^{-}\B_{\xi_{1}}^{+})=1.$$
If instead $2<j<p$, then 
$\ve_{l}\in \{0, +1\}$, for  all $ 2\leq l\leq j-1$, the following two cases appear:
\begin{enumerate}
\item There exists $2 \leq l \leq j-1$ such that $\ve_{l}=0$. Here, we define 
$$l'=\min\{l: 2\leq l\leq j-1, \ve_{l}=0\},$$
and
$$\vp(\B_{\bm\ve}^{\bm\xi})=\vp(\B_{\xi_{p}}^{\ve_{p}}\B_{\xi_{p-1}}^{\ve_{p-1}}\cdots \B_{\xi_{j+1}}^{\ve_{j+1}} \B_{\xi_{j}}^{-} \B_{\xi_{j-1}}^{\ve_{j-1}}\cdots \B_{\xi_{l'+1}}^{\ve_{l'+1}}\B_{\xi_{l'}}^{\circ}\B_{\xi_{l'-1}}^{+}\cdots B^{+}_{\xi_{2}}B^{+}_{\xi_{1}}\Omega)\neq 0.$$ 
Therefore, $$\B_{\xi_{j}}^{-} \B_{\xi_{j-1}}^{\ve_{j-1}}\cdots \B_{\xi_{l'+1}}^{\ve_{l'+1}}\B_{\xi_{l'}}^{\circ}\B_{\xi_{l'-1}}^{+}\cdots B^{+}_{\xi_{2}}B^{+}_{\xi_{1}}\Omega\neq 0.$$ \\
By \eqref{cr1} and \eqref{cr3}, one has 
$$\xi_{l'}=\xi_{l'-1}\succ \xi_{l'-2}\cdots \succ \xi_{2}\succ \xi_{1}$$
and
 \begin{align*}
\B_{\xi_{j}}^{-}\B_{\xi_{j-1}}^{\ve_{j-1}}\cdots \B_{\xi_{l'}}^{\circ}\B_{\xi_{l'-1}}^{+}\cdots \B_{\xi_{2}}^{+}\B_{\xi_{1}}^{+}\Omega &=\B_{\xi_{j}}^{-}\B_{\xi_{j-1}}^{\ve_{j-1}}\cdots \B_{\xi_{l'+1}}^{\ve_{l'+1}}\B_{\xi_{l'}}^{\circ}(g_{\xi_{l'-1}}\otimes \cdots\otimes g_{\xi_{1}})\\
  &=\B_{\xi_{j}}^{-}\B_{\xi_{j-1}}^{\ve_{j-1}}\cdots \B_{\xi_{l'+1}}^{\ve_{l'+1}}(g_{\xi_{l'-1}}\otimes \cdots\otimes g_{\xi_{1}})\\
  &=\B_{\xi_{j}}^{-}\B_{\xi_{j-1}}^{\ve_{j-1}}\cdots \B_{\xi_{l'+1}}^{\ve_{l'+1}}\B_{\xi_{l'-1}}^{+}\cdots \B_{\xi_{2}}^{+}\B_{\xi_{1}}^{+}\Omega\neq 0.
 \end{align*}
Moreover, the action of $\B_{\bm\xi}^{\bm\ve}$ on $\Omega$ is the same as the action of $\B_{\bm\xi'}^{\bm\ve'}$ on $\Omega$, where $\bm\ve'=(\ve_{p}, \ldots, \ve_{l'+1}, \ve_{l'-1}, \ldots,  \ve_{1})\in \D_{p-1}, \bm\xi'=(\xi_{p}, \ldots, \xi_{l'+1}, \xi_{l'-1}, \ldots, \xi_{1})\in [0, \rho]_{{\I}}^{p-1}$ and $\bm\xi'\sim \pi' \in \NC_{2}^{1, i}(p-1)$. We can now use the induction hypothesis to obtain that  $\bm\xi'$ establishes bm-order on $\pi'=\pi\setminus\{l'\} \in \NC_{2}^{1, i}(p-1) $. Moreover, since the block singleton $\{l'\}$ is a direct successor of the  pair block $\{l'-1, T(l'-1)\}$ and  $\xi_{l'}=\xi_{l'-1} \succ \cdots \succ \xi_{1}$, then 
 $\bm\xi\trianglelefteq\pi\in \NC_{2}^{1, i}(p).$
 \item For all $2\leq l \leq j-1, \ve_{l}=+1$. In this case, $j=T(j-1)$. Hence,
$$\vp(\B_{\bm\xi}^{\bm\ve})\neq 0\Longrightarrow \B_{\xi_{j}}^{-}\B_{\xi_{j-1}}^{+}\cdots \B_{\xi_{2}}^{+}\B_{\xi_{2}}^{+}\Omega\neq 0$$ and
$$\xi_{j}=\xi_{j-1}\succ \xi_{j-2}\succ \cdots \succ \xi_{1}.$$
Consequently, 
\begin{align*}
\B_{\xi_{j}}^{-}\B_{\xi_{j-1}}^{+}\cdots \B_{\xi_{1}}^{+}\Omega &=\B_{\xi_{j}}^{-}\B_{\xi_{j-1}}^{+}(g_{\xi_{j-2}}\otimes \cdots \otimes g_{\xi_{1}})\\
&=g_{\xi_{j-2}}\otimes g_{\xi_{j-3}}\otimes \cdots \otimes g_{\xi_{1}}\\
&=\B_{\xi_{j-2}}^{+}\B_{\xi_{j-3}}^{+}\cdots \B_{\xi_{1}}^{+}\Omega,
\end{align*}
and $$\vp(\B_{\bm\xi}^{\bm\ve})=\vp(\B_{\xi_{p}}^{-}\cdots \B_{\xi_{j+1}}^{\ve_{j+1}}\B_{\xi_{j-2}}^{+}\cdots \B_{\xi_{1}}^{+})=\vp(\B_{\bm\xi'}^{\bm\ve'})\neq 0,$$ 
where 
$$\bm\xi':=(\xi_{p}, \ldots, \xi_{j+1}, \xi_{j-2}, \ldots, \xi_{1})\sim \pi'\in\NC_{2}^{1, i}(p-2)$$
and
$$\bm\ve':=(\ve_{p}, \ldots, \ve_{j+1}, \ve_{j-2}, \ldots, \ve_{1})\in\D_{p-2}, \text{ where } \ve_1=\cdots=\ve_{j-2}=+1.$$
The induction assumption applied to the sequences $\bm\ve'$ and $\bm\xi'$, implies
$$\bm\xi'\trianglelefteq \pi'=\pi\setminus\{j-1, j\}\in\NC_{2}^{1, i}(p-2).$$
Since $\xi_{j}=\xi_{j-1}\succ \xi_{j-2}\cdots \succ \xi_{1}$, we have $\bm\xi\trianglelefteq \pi\in \NC_{2}^{1, i}(p)$.
\end{enumerate}
\end{enumerate}
\eproof
\begin{rem}\label{phi1}
The above arguments show that, for a sequence $\bm\xi=(\xi_{p}, \ldots, \xi_{1})\trianglelefteq \pi\in\NC_{2}^{1, i}(p)$, one has $\vp(\B_{\bm\xi}^{\bm\ve})=1$.
\end{rem}
By the previous results, the  proof of Theorem \ref{mthm} reduces to ascertain that for each $p\in {\N}$ the following limit exists 
$$\lim_{\rho\xrightarrow[]{\Pi} \infty} \frac{1}{\vo(\rho)^{\frac{p}{2}}}\sum_{\pi\in\NC_{2}^{1, i}(p)}\hspace{0.1cm}\sum_{\bm \xi\in \mathrm{bmo}(\pi, \rho)}\vp(\B_{\xi_{p}}^{\ve_{p}}\cdots \B_{\xi_{1}}^{\ve_{1}})\lambda_{\rho}^{s(\pi)}=\sum_{\pi \in \NC_{2}^{1, i}(p)} V(\tilde \pi)\lambda^{s(\pi)},$$
where $V(\tilde \pi)$ satisfies the recurrence relation  given in \eqref{rfv}. Remark \ref{phi1} yields that 
$$\lim_{\rho \xrightarrow[]{\Pi}\infty}\frac{1}{\vo(\rho)^{\frac{p}{2}}}\sum_{\pi\in \NC_{2}^{1, i}(p)}\hspace{0.1cm}\sum_{\bm \xi\in {\bmo}(\pi, \rho)}\vp(\B_{\xi_{p}}^{\ve_{p}}\cdots \B_{\xi_{1}}^{\ve_{1}})\lambda_{\rho}^{s(\pi)}=\sum_{\pi\in \NC_{2}^{1, i}(p)}\hspace{0.1cm}\lim_{\rho \xrightarrow []{\Pi}\infty}\frac{|\mathrm{bmo}(\pi, \rho)|}{\vo(\rho)^{b(\tilde \pi)}}\lambda^{s( \pi)},$$
where $b(\tilde \pi)=\frac{p-s(\pi)}{2}$ is the number of blocks of the reduced partition $\tilde \pi\in \NC_{2}(p-s(\pi))$ obtained by $\pi\in\NC_{2}^{1, i}(p)$ after extracting the singleton blocks. Theorem \ref{mthm} thus follows from the following.
\begin{prop}\label{limv}
Let $\pi\in \NC_{2}^{1, i}(p)$ be a noncrossing partition with pair or  inner singleton blocks. Then for each positive symmetric cone $\Pi_d$
$$\lim_{\rho\xrightarrow[]{\Pi} \infty}\frac{|\mathrm{bmo}(\pi, \rho)|}{\vo(\rho)^{b(\tilde \pi)}}= V (\tilde \pi),$$
where the function $V(\tilde \pi):=V_{\Pi_d}(\tilde \pi)$ depends on the cone $\Pi_d$ and can be recursively expressed by the volume characteristic sequence $(\gamma_{n}(\Pi_d))_{n\geq 0}$ as in \eqref{rfv}. 
\end{prop}
\bproof
By Theorem \ref{thmp}, we need to prove that
\begin{equation}\label{bmoeq}
|\mathrm {bmo}(\pi, \rho)|=|\mathtt{BMO}(\tilde \pi, \rho)|.
\end{equation}
Let us consider the sequence $\bm\xi$ with the adapted partition $\pi=(B_{1}, \ldots, B_{\frac{p+s(\pi)}{2}})\in \NC_{2}^{1, i}(p; \frac{p+s(\pi)}{2})$, where $B_{1}, \ldots, B_{\frac{p+s(\pi)}{2}}$ are the $\frac{p+s(\pi)}{2}=k$ blocks of $\pi$, and denote by $\mu:=(\mu_k, \ldots, \mu_1)$ the related label sequence. Since $[0, \rho]_{\I}^{p}\ni\bm\xi\trianglelefteq \pi$, then $[0, \rho]_{\I}^{k}\ni\mu\trianglelefteq_{\mathtt{L}}\pi$, and by Remark {\ref{equalsets}}, we obtain
$$|\mathrm{bmo}(\pi, \rho)|=|\mathrm{BMO}(\pi, \rho)|.$$
Moreover, the labels of the direct successors (singleton blocks) are the same as those of its direct predecessors (pair blocks). Then
$$|\mathtt{BMO}(\tilde \pi, \rho)|=|\mathrm{BMO}(\pi, \rho)|=|\mathrm{bmo}(\pi, \rho)|.$$ 
\eproof
\section*{Appendix}

In this section, we present the results for the single operator $S(\lambda):=A_{\xi}^{+}+A_{\xi}^{-}+\lambda A_{\xi}^{\circ}$, for $\lambda>0$ and $\xi\in{\I}$. In particular, the $p$-th moment $a_{p}(\lambda):=\vp((S(\lambda))^{p})$ and the associated probability measure $\nu_\lambda$ are given. It is worthwhile to mention that they are the same as those obtained by Muraki \cite{Mur0}.
\begin{prop}\label{propap}
 For any positive integer $p\geq 2$, one has 
 \begin{equation}\label{ap}
  a_{p}(\lambda)=\sum_{\pi\in \NC_{2, o}^{1, i}(p)}\lambda^{s(\pi)},
 \end{equation}
where $\NC_{2, o}^{1, i}(p)\subset\NC_{2}^{1}(p)$ is the set of all noncrossing  partitions with pair or singleton blocks such that the pair blocks must be outer and the singletons must be inner.\\
Furthermore, the following recursive formula holds
 \begin{equation}\label{ra}
  a_{0}(\lambda)=1, a_{1}(\lambda)=0, a_{p}(\lambda)=\lambda a_{p-1}(\lambda)+a_{p-2}(\lambda), \quad \text{ for }\quad  p\geq 2.
 \end{equation}
\end{prop}

\begin{exmp}
For $p=0, 1, \ldots, 6$ we can use the recursive formula \eqref{ra}, to obtain the following moment sequence:
$$a_{p}(\lambda)=\{1, 0, 1, \lambda, \lambda^2+1, \lambda^3+2\lambda, \lambda^4+3\lambda^2+1\}.$$
\end{exmp}

The moment generating function and the Cauchy transform for the moments of $S(\lambda)=A_{\xi}^{+}+A_{\xi}^{-}+\lambda A_{\xi}^{\circ}$ are given, respectively, by 
 \begin{equation}\label{mgf}
  \mathrm{M}_{\lambda}(x)=\frac{1-\lambda x}{1-\lambda x- x^{2}},
 \end{equation}
 and 
\begin{equation}\label{ctr}
 \mathrm{G}_{\lambda}(x)=\frac{x-1}{x^{2}-\lambda x-1}.
\end{equation}
Then, by the Stieltjes inverse formula, the probability distribution $\nu_{\lambda}$ of $ S(\lambda)$ under the vacuum state $\vp$ is given as follows
$$\nu_{\lambda}=p_{1}\delta_{x_{1}}+p_{2}\delta_{x_{2}},$$ 
where $x_{1}=\frac{\lambda}{2}+\frac{\sqrt{\lambda^{2}+4}}{2}$, $x_{2}=\frac{\lambda}{2}-\frac{\sqrt{\lambda^{2}+4}}{2}$, $p_{1}=\frac{1}{2}-\frac{\lambda}{2\sqrt{\lambda^{2}+4}}$, $p_{2}=\frac{1}{2}+\frac{\lambda}{2\sqrt{\lambda^{2}+4}}$ and  $\delta_{x}$ denotes the Dirac measure at a point $x$.

Furthermore, according to Corollary \ref{bmc}, one can show that the operators $\{S(\lambda)\}_{\xi}$ are bm-independent with respect to $\vp$ in $(\mathcal{C}, \vp)$, and then the operators $\{S(\lambda)\}_{\xi}$ can be considered as \emph{bm-Bernoulli} random variables.
\begin{rem}
For $\lambda=1$ and $p\geq 1$, the moment sequence $a_{p}(1)$ are the shifted Fibonacci numbers, given by
 $$a_{p}(1)=F(p-1).$$
\end{rem}
\vspace{1cm}

\section*{Acknowledgements}

Research partially supported by the National Agency for Academic Exchange (NAWA) POLONIUM project PPN/BIL/2018/1/00197/U/00021 and by the Polish National Science Center (NCN) grant 2016/21/B/ST1/00628.

\end{document}